\documentclass[11pt, reqno]{amsart}
\usepackage[margin=1in]{geometry}
\usepackage{graphicx}
\numberwithin{equation}{section}
\usepackage[T1]{fontenc}
\usepackage{color}
\usepackage{amsmath}
\usepackage{enumitem} 
\usepackage{amssymb}
\usepackage{mathtools}
\usepackage{commath}
\usepackage{hyperref}
\usepackage[numbers,sort]{natbib}

\newcommand{\iid}{\text{i.i.d}}
\newcommand{\Var}{\mathrm{Var}}

\newcommand{\R}{\mathbb{R}}
\newcommand{\one}{{\bf 1}}
\def\P{\mathbb{P}}
\def\E{\mathbb{E}}

\newcommand{\N}{\mathbb{N}}

\newcommand{\X}{\mathcal{X}}
\newcommand{\Y}{\mathcal{Y}}
\newcommand{\K}{\mathcal{K}}

\newcommand{\ind}[1]{\mathbf{1}\big\{#1\big\}} 
\newcommand{\y}{\mathbf{y}}
\newcommand{\x}{\mathbf{x}}

\newcommand{\diam}{\mathrm{diam}}
\newcommand{\card}[1]{|#1|}
\DeclarePairedDelimiterX{\inp}[2]{\langle}{\rangle}{#1, #2}

\makeatletter
\newcommand{\leqnos}{\tagsleft@true\let\veqno\@@leqno}

\newcommand{\Xn}{\mathcal{X}_n}

\newcommand{\F}{\mathcal F}
\newcommand{\stan}{R_{q_m} +a(R_{q_m})\rho}

\makeatletter
\@namedef{subjclassname@2020}{%
  \textup{2020} Mathematics Subject Classification}
\makeatother

\theoremstyle{plain}
\newtheorem{theorem}{Theorem}[section]
\newtheorem{lemma}[theorem]{Lemma}
\newtheorem{proposition}[theorem]{Proposition}

\theoremstyle{definition}
\newtheorem{definition}[theorem]{Definition}

\newtheorem{example}[theorem]{Example}  

\title[FSLLNs for Euler characteristic processes of extreme sample clouds]{Functional strong laws of large numbers for Euler characteristic processes of extreme sample clouds}
\author{Andrew M. Thomas \and Takashi Owada}
\date{}

\address{Department of Statistics\\
Purdue University \\
IN, 47907, USA}
\email{thoma186@purdue.edu \\ owada@purdue.edu}


\begin{document}

\maketitle

\begin{abstract}
To recover the topology of a manifold in the presence of heavy tailed or exponentially decaying noise, one must understand the behavior of geometric complexes whose points lie in the tail of these noise distributions. This study advances this line of inquiry, and demonstrates functional strong laws of large numbers for the Euler characteristic process of random geometric complexes formed by random points outside of an expanding ball in $\mathbb{R}^d$. When the points are drawn from a heavy tailed distribution with a regularly varying tail, the Euler characteristic process grows at a regularly varying rate, and the scaled process converges uniformly and almost surely to a smooth function. When the points are drawn from a distribution with an exponentially decaying tail, the Euler characteristic process grows logarithmically, and the scaled process converges to another smooth function in the same sense. All of the limit theorems take place when the points inside the expanding ball are densely distributed, so that the simplex counts outside of the ball of all dimensions contribute to the Euler characteristic process. 
\end{abstract} 

\maketitle

\section{Introduction}

\subsection{Heavy tailed noise and annuli structure of homological elements}  \label{sec:noise.annuli}

\begin{figure}[t]
\centering
\includegraphics[width = 4in]{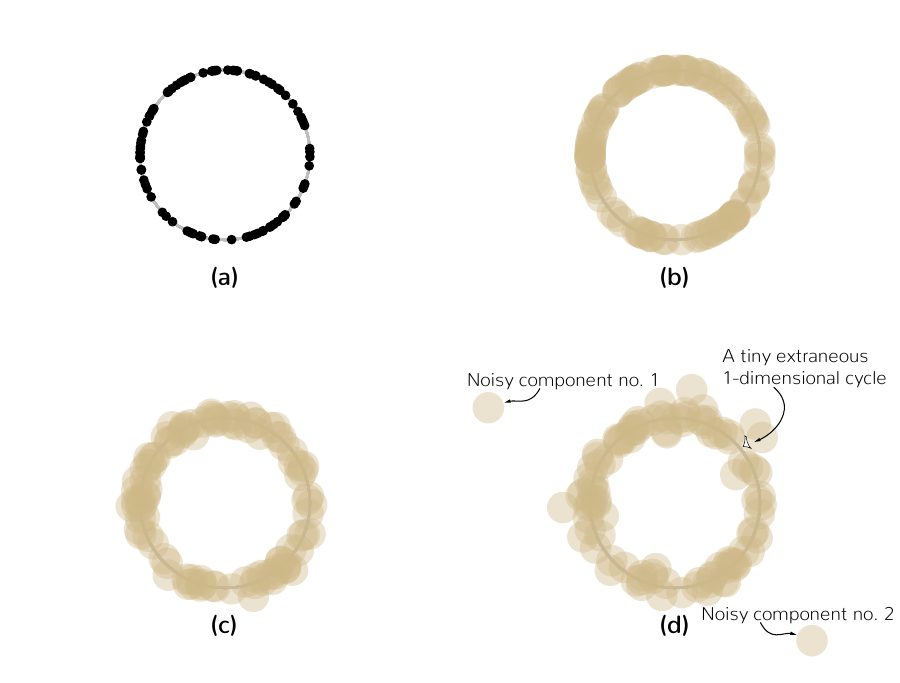}
\caption{(a) A random sample $\X_{100} = \{ X_1,\dots,X_{100} \}$ of 100 points, uniformly distributed on the unit circle $S^1$. (b) The union of balls of radius 0.2 around each point in $\X_{100}$. (c) The union of balls when $\X_{100}$ is perturbed by Gaussian noise. (c) The union of balls when $\X_{100}$ is perturbed by Cauchy noise. Noise is applied in the same manner as \cite{nys2011}. Topology is recovered by the union of the balls in cases (b) and (c), but that is not the case if the heavy tailed Cauchy noise is added.}
\label{f:perils_ht_noise}
\end{figure}

To recover the topology of a manifold using point cloud data, one needs to have a strong understanding of how the points are perturbed from the manifold. In \cite{nys2008}, given a ``nice'' manifold, it was shown that one can recover the topology of the manifold by a sufficiently dense random sampling of points if the noise is bounded. In \cite{nys2011} it was shown that the recovery is still possible by a sufficiently dense random sampling of points if the noise is standard multivariate Gaussian and the variance is bounded by a function of the reach and dimension of the manifold. 
However, if the points on the manifold are perturbed by heavy tailed noise, the recovery of topology will be severely impacted because of extraneous homological elements generated by this noise. In Figure~\ref{f:perils_ht_noise}, we wish to recover the topology of the circle $S^1$ from the union of balls around a random point cloud. In Figure \ref{f:perils_ht_noise} (b) and (c), the union of balls recover the essential shape of $S^1$, as the size of noise in these cases is sufficiently small. However, in Figure \ref{f:perils_ht_noise} (d) the noise added to points in $S^1$ has a heavy tailed Cauchy distribution. Consequently, three extraneous shape elements appear --- two distinct components and a tiny one-dimensional cycle. This phenomenon in case (d) raises the question of how the shape of these elements away from the center of $S^1$ may behave in general; this is roughly the idea of what is called \emph{topological crackle}.

From a more analytic viewpoint, topological crackle is understood as a layered annuli structure of homological elements of different orders. 
To make our implication more clear, we consider the power-law density, 
\begin{equation}  \label{e:power.law.pdf.intro}
f(x) = \frac{C}{1+\|x\|^\alpha}, \ \ \ x\in \R^d, 
\end{equation}
for some $\alpha >d$ and a normalizing constant $C>0$. Suppose a random point cloud $\X_n=\{ X_1,\dots,X_n \} \subset \R^d$, $d\ge 2$, is drawn from this density. Let Ann$(K,L)$ be a closed annulus with inner radius $K$ and outer radius $L$, and $B(x,t):=\big\{ y\in\R^d: \|y-x\|<t \big\}$ be an open ball of radius $t$ around $x\in \R^d$ (here $\|\cdot\|$ denotes the Euclidean norm). Then, one can divide $\R^d$ in a way that 
\begin{equation}  \label{e:annuli}
\R^d = \bigcup_{k=0}^{d+1} \text{Ann} (R_{k,n}, R_{k-1,n}), 
\end{equation}
where
\begin{equation}  \label{e:annuli.radii}
R_{k,n} = \begin{cases}
\infty & \ k=-1, \\
(Cn)^{\frac{1}{\alpha-d}} & \ k=0, \\
(Cn)^{\frac{1}{\alpha-d/(k+2)}} & \ k\in \{ 1,\dots,d-1 \}, \\
(Cn)^{1/\alpha} & \ k=d, \\
0 & \ k=d+1, 
\end{cases}
\end{equation}
so that 
$$
 R_{d,n} \ll R_{d-1,n} \ll \cdots \ll R_{1,n} \ll R_{0,n}, \ \ \text{as } n\to\infty. 
$$
From previous studies \citep{crackle, owada_crackle, owada2018}, it is known that the union of balls, 
$$
U(t) := \bigcup_{X\in \X_n} B(X,t), \ \ \ t\ge 0, 
$$
asymptotically generates homological elements (i.e., components and cycles as seen in Figure \ref{f:perils_ht_noise} (d)) in the following way: we have, as $n\to\infty$, 
\vspace{3pt}

\begin{itemize}
\item Inside Ann$(R_{0,n}, \infty) = B(0,R_{0,n})^c$ there are finitely many distinct components, but none of the cycles of dimensions $1,2,\dots,d-1$. 
\vspace{3pt}

\item Inside Ann$(R_{1,n}, R_{0,n})$ there are infinitely many distinct components and finitely many one-dimensional cycles, but none of the cycles of dimensions $2,3,\dots,d-1$. 
\end{itemize}
\vspace{3pt}
In general, for every $k\in \{ 2,\dots,d-1 \}$, 
\vspace{3pt}

\begin{itemize}
\item Inside Ann$(R_{k,n}, R_{k-1,n})$ there are infinitely many distinct components and cycles of dimensions $1,\dots,k-1$, and finitely many $k$-dimensional cycles, but none of the cycles of dimensions $k+1,\dots,d-1$, 
\end{itemize}
\vspace{3pt}
and finally,
\vspace{3pt}
 
\begin{itemize}
\item Inside Ann$(R_{d,n}, R_{d-1,n})$ there are infinitely many distinct components and cycles \\ of \emph{all} dimensions $1,\dots,d-1$. 
\end{itemize}
\vspace{3pt}

In the literature \citep{owada2017}, the innermost ball $B(0,R_{d,n})$ is referred to as a \emph{weak core}. Inside the weak core, random points are densely scattered and the homology of the union of unit balls around them is nearly trivial as $n\to\infty$, i.e., the union has nearly no cycles of all dimensions $1,\dots,d-1$. With this layered structure in mind, topological crackle is formally defined as non-trivial homological elements (i.e., components and cycles) outside of a weak core. 

Many of the existing studies in classical extreme value theory (EVT) have focused on the behavior of a random point cloud in the outermost annulus Ann$(R_{0,n}, \infty)$, or equivalently outside of $B(0,R_{0,n})$. For instance, it is well known that the total number of distinct components outside of $B(0,R_{0,n})$ converges weakly to a Poisson distribution as $n\to\infty$ \citep{resnick:1987, htresnick, embrechts:kluppelberg:mikosch:1997}. 
The objective of this paper is to go beyond the studies on the spatial distribution of components and investigate more complicated and higher-dimensional topological features outside of a weak core. 

\subsection{Topological crackle and Euler characteristic} 
After the pioneering paper of \cite{crackle}, the stochastic properties of crackle phenomena have been investigated mostly via the behavior of \emph{Betti numbers}. Loosely speaking, the $k$th Betti number counts the number of $k$-dimensional cycles which can be interpreted as the boundary of a $(k+1)$-dimensional body. In the related literature, 
\cite{owada_crackle} studied the case in which the $k$th Betti number of $U(t)$ outside of $B(0,R_{k,n})$ converges weakly to a Poisson distribution as $n\to\infty$. Moreover, \cite{owada2018} established the central limit theorem for the $k$th Betti numbers, in the case that infinitely many $k$-dimensional  cycles appear outside of an expanding ball. Additionally, \cite{owada:2019} gave a rigorous description of the limiting Betti numbers when the random points are generated by a classical moving average process, and \cite{owada2020} discussed the weak convergence of a standard graphical representation of cycles. 

In contrast to these previous papers, the objective of this paper is to examine the crackle phenomena from the viewpoint of the \emph{Euler characteristic} of a geometric (simplicial) complex. 
Among many varieties of geometric complexes \citep[see][]{ghrist:2014}, the  \emph{Vietoris-Rips complex} and the \emph{\v{C}ech complex} are specific examples that deserve our attention. 
\begin{definition}  \label{def:VR}
Given a point set $\X=\{x_1,\dots,x_n\} \subset \R^d$ and a positive number $t>0$, the Vietoris-Rips complex $\mathcal R(\X,t)$ is defined as follows. 
\begin{itemize}
\item The $0$-simplices are the points in $\X$. 
\item A $k$-simplex $\sigma=[x_{i_0}, \dots,x_{i_k}]$ is in $\mathcal R(\X,t)$ if $\bar{B}(x_{i_p},t/2) \cap \bar{B}(x_{i_q}, t/2)\neq \emptyset$ for every $0\le p<q\le k$, where $\bar{B}(x,t)$ is the closure of $B(x,t)$. 
\end{itemize}
\end{definition}
\begin{definition}  \label{def:Cech}
Given the same  $\X$ and $t>0$, the \v{C}ech complex $\check C(\X,t)$ is defined as follows. 
\begin{itemize}
\item The $0$-simplices are the points in $\X$. 
\item A $k$-simplex $\sigma=[x_{i_0}, \dots,x_{i_k}]$ is in $\check C(\X,t)$ if a family of balls $\big\{  \bar{B}(x_{i_p},t/2), \, p=0,\dots,k \big\}$ has a non-empty intersection. 
\end{itemize}
\end{definition}



In this paper, we examine two distinct scenarios of noise distributions that experience topological crackle: one where the distribution has a regularly-varying tail and another where the distribution has an exponentially decaying tail. We define a geometric complex that generalizes both $\mathcal R(\X,t)$ and $\check C(\X,t)$ above, and then establish the functional strong law of large numbers (FSLLN) for the corresponding Euler characteristic process for each of the distributional contexts. One of the primary benefits of working with the Euler characteristic comes from the fact that it can be expressed as an alternating sum of Betti numbers of all dimensions \citep{edelsbrunner}. As a consequence, the Euler characteristic provides a limit theorem containing information on cycles of all different dimensions, whereas Betti numbers can only provide separate limit theorems for cycles of a particular dimension in each individual annulus region of \eqref{e:annuli}. Such global features of the Euler characteristic help to capture the spatial distribution of cycles in the union of annuli $\bigcup_{k=0}^d \text{Ann}(R_{k,n}, R_{k-1,n}) = \R^d \setminus B(0,R_{d,n})$, even though the nature of the distribution of cycles differs from region to region.

In conjunction with the recent development of Topological Data Analysis (TDA), the literature dealing with the asymptotics of the Euler characteristic of random geometric complexes has flourished \citep{thomas:owada:2020, krebs2020, bobrowski:adler:2014, bobrowski:mukherjee:2015, hug:last:schulte:2016}. However, none of these studies have paid sufficient attention to the topology of the tail of a probability distribution. In the context of topological crackle as in Figure \ref{f:perils_ht_noise}, ascertaining the topology of noise is an important step in determining how to process the signal of the manifold. In the light of the results in this paper, one can examine how the homology (i.e., components and cycles) of extreme-valued noise distributions evolve. This can be attained by viewing the Euler characteristic as a stochastic process, in which the parameter governing the formation of simplices is taken to be the ``time'' parameter. The resulting process then relates strongly to persistent homology. Persistent homology is a topological and algebraic structure that tracks the creation and destruction of cycles in different dimensions. It is one of the most widely used and robust tools in the TDA toolbox --- see \cite{adler2010pers} or \cite{ghrist:2008} for an introduction to persistent homology and \cite{edelsbrunner} for a more thorough treatment. In particular, Examples \ref{ex:chi_np} and \ref{ex:exp.density} below provide the FSLLNs in the different scenarios of noise distributions for the integrated Euler characteristic process. This process can be viewed as the Euler characteristic of a persistence barcode, which is a well known graphical descriptor of persistent homology \citep{ghrist:2008, carlsson:2009}. 

\subsection{Organization of the paper} 
The remainder of this paper is structured as follows. Section \ref{s:prelim} provides a discussion of the background material necessary for this paper. The paper then proceeds to the heavy tailed setup and presents the FSLLN for the Euler characteristic in Section \ref{s:heavy}. The paper continues with a discussion of the intricacies of the exponentially decaying tail case along with the corresponding FSLLN in Section~\ref{s:lite}. The proofs of the main results for both setups are deferred to Section \ref{s:proof}. From a technical point of view, the studies most relevant to this paper are \cite{goel2018} and  \cite{thomas:owada:2020}, in which the authors established strong laws of large numbers for topological invariants --- such as Betti numbers and the Euler characteristic --- in the non-extreme value theoretic setup. In particular, these studies revealed that if the topological invariants are scaled proportionally to the sample size, they converge almost surely to a finite and positive constant. Owing to this fact, the main machinery in their proofs is a direct application of the Borel-Cantelli lemma, together with the calculation of lower-order moments. On the contrary, the main challenge in this paper is that the scaling sequence of the Euler characteristic may grow very slowly (e.g., logarithmically), in which case, a direct application of the Borel-Cantelli lemma does not work. To overcome this difficulty, we need to identify suitable subsequential upper and lower bounds of the Euler characteristic to which one can apply the Borel-Cantelli lemma. This is a standard technique in the theory of random geometric graphs --- see Chapter 3 of the monograph of \cite{penr}. It is possible to extend these arguments to our higher-dimensional setup since the geometric complexes such as those in Definitions \ref{def:VR} and \ref{def:Cech} are higher-dimensional analogues of a geometric graph.

As a final remark, we point out that the other types of limit theorems for the Euler characteristic still remain as a future topic. For example, it seems feasible to establish a (functional) central limit theorem for the Euler characteristic via Stein's method for normal approximation \citep[see, e.g., Theorem 2.4 in][]{penr}. Indeed, \cite{owada2018} already derived the central limit theorem for the Betti numbers by means of the aforementioned normal approximation technique. We anticipate that the same approach is possible for our Euler characteristic.  A similar line of research in the non-extreme value theoretic setup can be found in \cite{thomas:owada:2020} and \cite{krebs2020}.

\section{Preliminaries} \label{s:prelim}

The point cloud of interest in this study is the sample $\Xn := \{X_1, \dots, X_n\}$ of $n$ $\iid$ random points in $\R^d$, $d \geq 2$ with spherically symmetric density $f$. Spherical symmetry of $f$ is far from necessary; the results in this paper could be extended to densities with ellipsoidal level sets fairly easily. Because of the imposed spherical symmetry, we can define $f(r):=f(r\theta)$ for all $r\ge 0$ and $\theta \in S^{d-1}$.
Denote $\lambda$ to be Lebesgue measure on $\R^d$ and $S^{d-1} := \{x \in \R^d: \norm{x} = 1\}$. Let us here define the spherical measure
\[
\nu_{d-1}(A) := d\cdot \lambda \big(\big\{ x\in B(0,1): x/\|x\|\in A \big\}\big)
\]
for Borel sets $A \subset S^{d-1}$. We denote $\omega_d := \lambda(B(0,1)) = 2\pi^{d/2}/ \big( d\Gamma(d/2) \big)$, and $s_{d-1} := \nu_{d-1}(S^{d-1}) = d\omega_d$.

Let $\F(\R^d)$ be the collection of all non-empty, finite subsets of $\R^d$. For $\X\in \F(\R^d)$, a \emph{simplicial complex} $\K(\X)$ is a collection of subsets of $\X$ such that if $\sigma \in \K(\X)$ and $\tau \subset \sigma$ then $\tau \in \K(\X)$. Evidently, the Vietoris-Rips complex $\mathcal R(\X,t)$ and the \v{C}ech complex $\check{C}(\X,t)$ satisfy this condition. 
We call $\sigma \in \K(\X)$ a \emph{k-simplex} if $\card{\sigma} = k+1$. 

Subsequently, let $h:\F(\R^d)\to \{ 0,1\}$ be an indicator function satisfying the following conditions. 
\begin{description}
\item[H1\label{H1}] $h(\X) \le h(\Y)$ for all $\Y\subset \X$. 
\item[H2\label{H2}] $h$ is translation invariant --- that is, for every $\X\in \F(\R^d)$ and $y\in\R^d$, we have $h(\X+y)=h(\X)$.
\item[H3\label{H3}] $h$ is locally determined --- that is, there exists $c>0$ so that $h(\X)=0$ whenever diam$(\X)>c$, where diam$(\X):= \max_{x,y\in \X}\|x-y \|$. 
\end{description}
By abusing notation slightly, for $\X=\{ x_1,\dots,x_m \}\in \F(\R^d)$, we write $h(\X) = h(x_1,\dots,x_m)$. Moreover, for $\X=\{ x_1,\dots,x_m \}$ and $a\in \R$, we write $a\X=\{ax_1, \dots, ax_m\}$. We then define a scaled version of $h$ by
$$
h_t(\X) := \begin{cases} h(t^{-1}\X), &  \ t>0,  \\
\ind{|\X| = 1}, & \ t=0, 
\end{cases}
$$
with the additional assumption that 
\begin{description}
\item[H4\label{H4}] $t \mapsto h_t(\X)$ is right continuous and non-decreasing for each $\X\in\F(\R^d)$.
\end{description}

Given such a scaled indicator $h_t$, we can construct the geometric (simplicial) complex
\begin{equation}  \label{e:simp.comp.general}
\K(\X, t) := \big\{ \Y\subset \X: h_t(\Y)=1  \big\}. 
\end{equation}
By virtue of \textbf{(H4)} above, \eqref{e:simp.comp.general} induces a \emph{filtration} of geometric complexes over a point set $\X$ --- that is, 
$$
\K(\X, s) \subset \K(\X, t) \ \ \text{for all } 0\le s \le t. 
$$
Note that if one takes 
\begin{equation*}  
h(\X) = \one \big\{ \text{diam}(\X) \le 1 \big\}, \ \ \X\in\F(\R^d), 
\end{equation*}
then \eqref{e:simp.comp.general} induces a \emph{Vietoris-Rips filtration}. Moreover, if we define 
$$
h(\X) = \one \Big\{\bigcap_{x \in \X} \bar{B}(x, 1/2) \neq \emptyset\Big\}, \ \ \X\in \F(\R^d), 
$$
then \eqref{e:simp.comp.general} induces a \emph{\v{C}ech filtration}. 

As mentioned in the Introduction, the objective of this paper is to study ``extreme-value'' behavior of random geometric complexes via the Euler characteristic. More concretely, with a non-random sequence $R_n\to \infty$, we study the filtration of geometric complexes
\begin{equation}  \label{e:geom.comp.EVT}
\K\big( \X_n \cap B(0,R_n)^c, t \big), \ \ \ t\ge 0, 
\end{equation}
which are distributed increasingly further from the origin as $n\to\infty$. We now define the Euler characteristic pertaining to \eqref{e:geom.comp.EVT} by 
\begin{equation}  \label{e:def.EC.proc}
\chi_n(t) := \sum_{k=0}^\infty (-1)^k S_{k,n}(t), \ \ \ t \ge 0, 
\end{equation}
where $S_{k,n}(t)$ denotes the $k$-simplex counts in the complex \eqref{e:geom.comp.EVT}. Namely, 
$$
S_{k,n}(t) := \sum_{\Y \subset \X_n, \, |\Y|=k+1} h_t(\Y)\, \one \big\{ \min_{y\in \Y} \| y\| \ge R_n \big\}. 
$$
Note that for every $n\in \N = \{1,2, \dots\}$, \eqref{e:def.EC.proc} is a finite sum as $S_{k,n}(t)\equiv 0$ for all $k\ge n$. 
Furthermore, \eqref{e:def.EC.proc} can be seen as a stochastic process in parameter $t$, with right continuous sample paths and left limits. In the following, we establish the FSLLN for the \emph{Euler characteristic process} $(\chi_n(t), \ t \geq 0)$ in the space $D[0,\infty)$ of right continuous functions on $[0,\infty)$ with left limits. In particular, we equip $D[0,\infty)$ with the uniform topology. 

With the notation in \eqref{e:def.EC.proc}, if we set $t=0$, 
\begin{equation}  \label{e:EC.proc.time0}
\chi_n(0) =\sum_{k=0}^\infty (-1)^k S_{k,n}(0) = \sum_{i=1}^n \one \big\{ \|X_i\|\ge R_n \big\}
\end{equation}
represents the number of points outside of an expanding ball $B(0,R_n)$. The asymptotics of \eqref{e:EC.proc.time0} can be treated in the standard framework of classical EVT \citep[see, e.g.,][]{resnick:1987,htresnick, embrechts:kluppelberg:mikosch:1997}. Unlike this special case, the Euler characteristic process \eqref{e:def.EC.proc} intrinsically involves higher-dimensional topological structures, which requires much more complicated machinery to analyze.

From the literature of topological crackle \citep[see][]{crackle, owada2017, owada2018, owada_crackle, owada2020}, it is known that the behavior of topological invariants significantly depends on the limit value of $nf(R_n)$. The present study focuses exclusively on the case when the limit of $nf(R_n)$ is a positive and finite constant --- that is, 
\begin{equation}  \label{e:weak.core.regime}
nf(R_n) \to \xi \ \ \text{as } n\to\infty \ \text{ for some } \xi \in (0,\infty).   
\end{equation}
As mentioned in the introduction, the ball $B(0,R_n)$ with $nf(R_n)\to 1$, $n\to\infty$, is called a weak core. In the special case when the density has a power-law tail as in \eqref{e:power.law.pdf.intro}, the radius of a weak core is equal to $(Cn)^{1/\alpha}$ (see \eqref{e:annuli.radii}). Therefore, if $R_n$ is determined by \eqref{e:weak.core.regime}, $B(0,R_n)$ coincides with the weak core, up to multiplicative constants. 
The configuration of points between the outside and inside of a weak core is very different. Inside of the weak core, the homology of the union of balls becomes almost trivial as $n \to \infty$, i.e., the random points are very densely distributed and nearly every cycle of every dimension becomes filled in \citep[see][]{crackle, owada2017}. Outside of the weak core, however, the random points are distributed more sparsely, though densely enough so that homology of all feasible dimensions becomes not only non-trivial, but abundant. As a consequence, by an appropriate scaling as a function of $R_n$ in \eqref{e:weak.core.regime}, the simplex counts of all dimensions in \eqref{e:def.EC.proc} will contribute to the limit. In contrast, if $R_n$ grows faster, such that $nf(R_n)\to0$ as $n\to\infty$, then even under an appropriate scaling, the Euler characteristic is dominated asymptotically by the 0-simplices, or the extremal points. In this case, the Euler characteristic simply counts points in $\X_n$ outside $B(0, R_n)$ as $n\to\infty$.

\section{Regularly varying tail case} \label{s:heavy}

In this section, we detail the large-sample behavior of \eqref{e:def.EC.proc} of an extreme sample cloud when the distribution of points has a heavy tail. Recall that $\X_n=\{X_1,\dots,X_n\}$ denotes a random sample in $\R^d$ with a spherically symmetric density $f$. We assume that there exists a tail index $\alpha >d$, such that 
\begin{equation}  \label{e:RV.tail}
\lim_{r \to \infty} \frac{f(rt\theta)}{f(r\theta)} = t^{-\alpha}, \ \ \forall t>0,
\end{equation}
for every (equivalently, some) $\theta\in S^{d-1}$. 

Before stating the main result, note that by \eqref{e:weak.core.regime}, one can typically take 
\begin{equation}  \label{e:explicit.Rn}
R_n = \xi^{-1/\alpha}\big(1/f\big)^{\leftarrow} (n),
\end{equation}
where $(1/f)^\leftarrow (x) = \inf\big\{  y: 1/f(y) \ge x \big\}$
is the (left continuous) inverse of $1/f$. 
Thus, $(R_n)_{n \geq 1}$ is a regularly varying sequence of exponent (tail index) $1/\alpha$.




\begin{theorem} \label{t:heavy}
Suppose that $f$ is a spherically symmetric density satisfying \eqref{e:RV.tail}. Assume that $nf(R_n) \to \xi$ as $n\to\infty$ for some $\xi\in (0,\infty)$. Then, the Euler characteristic process in \eqref{e:def.EC.proc} satisfies the following functional SLLN, i.e., as $n\to\infty$, 
\begin{equation}  \label{e:fslln.heavy}
\left( \frac{\chi_n(t)}{R_n^d}, \, t\ge 0  \right)\to \bigg( \sum_{k=0}^{\infty} (-1)^{k}s_{k}(t), \, t\ge0  \bigg),  \quad \mathrm{a.s. \ in \ } D[0, \infty).  
\end{equation}
where 
\begin{equation}  \label{e:limit.heavy}
s_k(t) := \frac{s_{d-1} \xi^{k+1}}{(k+1)!\big(\alpha(k+1)-d \big)} \int_{(\R^d)^k} h_t(0,y_1,\dots,y_k) \dif \y, \ \ \ t \ge 0, \ \ k\ge 1, 
\end{equation}
with $s_0(t)\equiv s_{d-1}\xi/(\alpha-d)$. In particular, the limit in \eqref{e:fslln.heavy} is convergent for all $t\ge 0$. 
\end{theorem}

The following example illustrates the uniform convergence that takes place in the above theorem.

\begin{example} \label{ex:chi_np}
Consider the power-law density defined by
$$
f(x) = \frac{2}{\pi \omega_d (1+ \norm{x}^{2d})}, \quad x \in \R^d. 
$$
Define $R_n := (2n/\pi\omega_d)^{1/(2d)}$, so that $nf(R_n) \to 1$. We consider the Vietoris-Rips complex induced by $h(\X) = \ind{\diam(\X) \leq 1/\sqrt{d}}$,  where $\diam$ is calculated here with respect to the $\ell^{\infty}$ norm. Then, it follows from Theorem \ref{t:heavy} that, as $n\to\infty$,
$$
\left( \sqrt{\frac{\pi\omega_d}{2n}}\chi_n(t), \, t\ge 0\right) \to \bigg(\sum_{k=0}^\infty (-1)^k s_k(t), \, t\ge 0\Big),  \ \ \text{a.s. in } D[0,\infty). 
$$
The limiting function above can be simplified as follows: 
\begin{align}
\sum_{k=0}^{\infty} (-1)^{k}s_{k}(t) &= s_{d-1} \sum_{k=0}^{\infty} \frac{(-1)^k(t/\sqrt{d})^{dk}}{ (k+1)! (2d(k+1) - d)} \int_{(\R^d)^k} h_{\sqrt{d}}(0, y_1, \dots, y_k) \dif{\y} \notag \\ 
&= \omega_d \sum_{k=0}^{\infty} \frac{(-1)^k (t/\sqrt{d})^{dk} (k+1)^d}{(k+1)! (2k+1)}. \label{e:ex_dens}
\end{align}
See Figures~\ref{fig:chi_np} and \ref{fig:multi.dim.EC} for actual plots of the limiting function in \eqref{e:ex_dens} for $d= 2,3,4,5$. 

\begin{figure}[t] 
\includegraphics[width=1.5in]{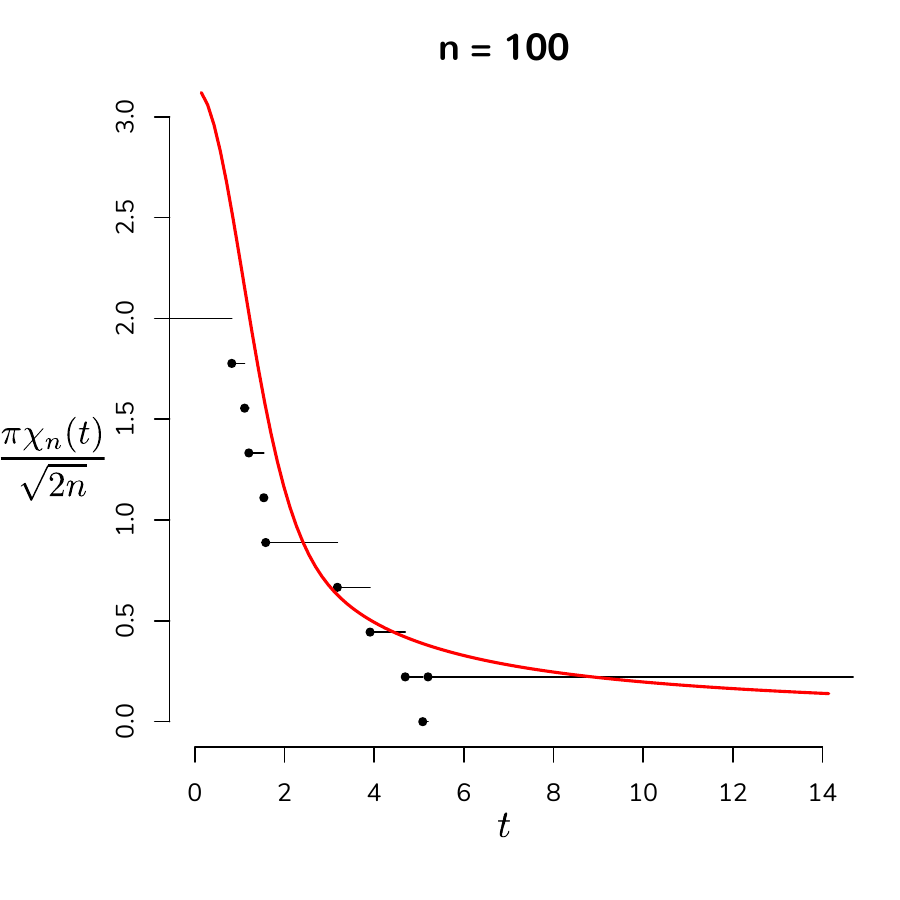} 
\includegraphics[width=1.5in]{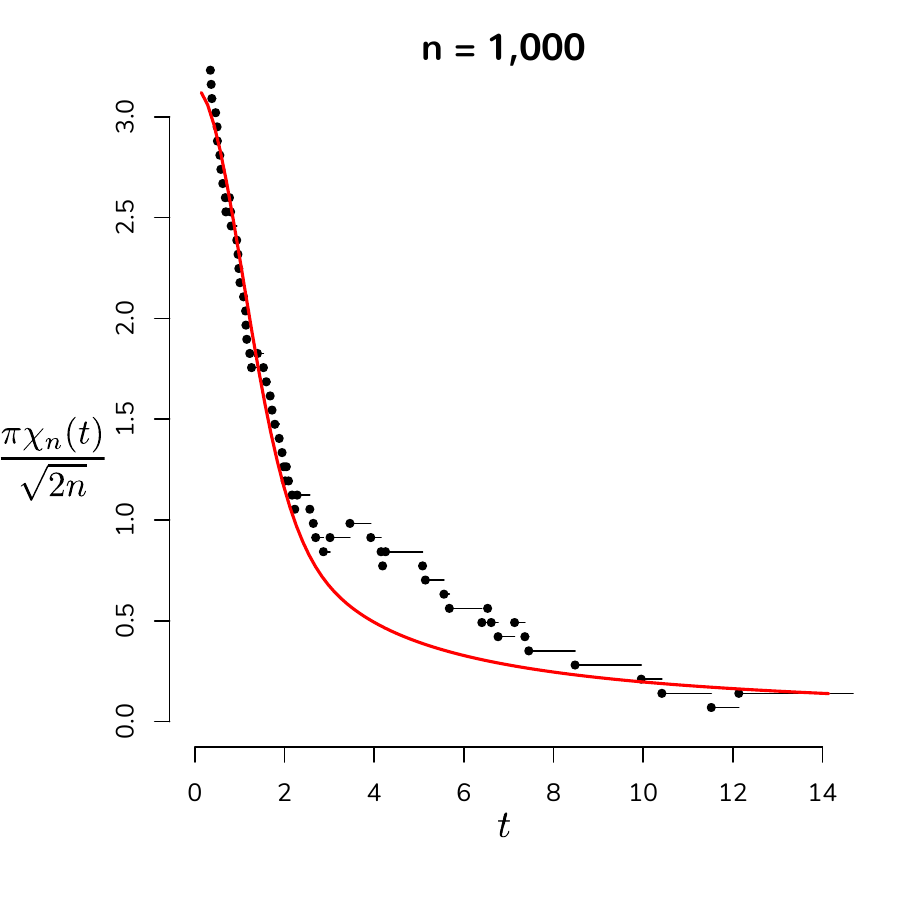}
\includegraphics[width=1.5in]{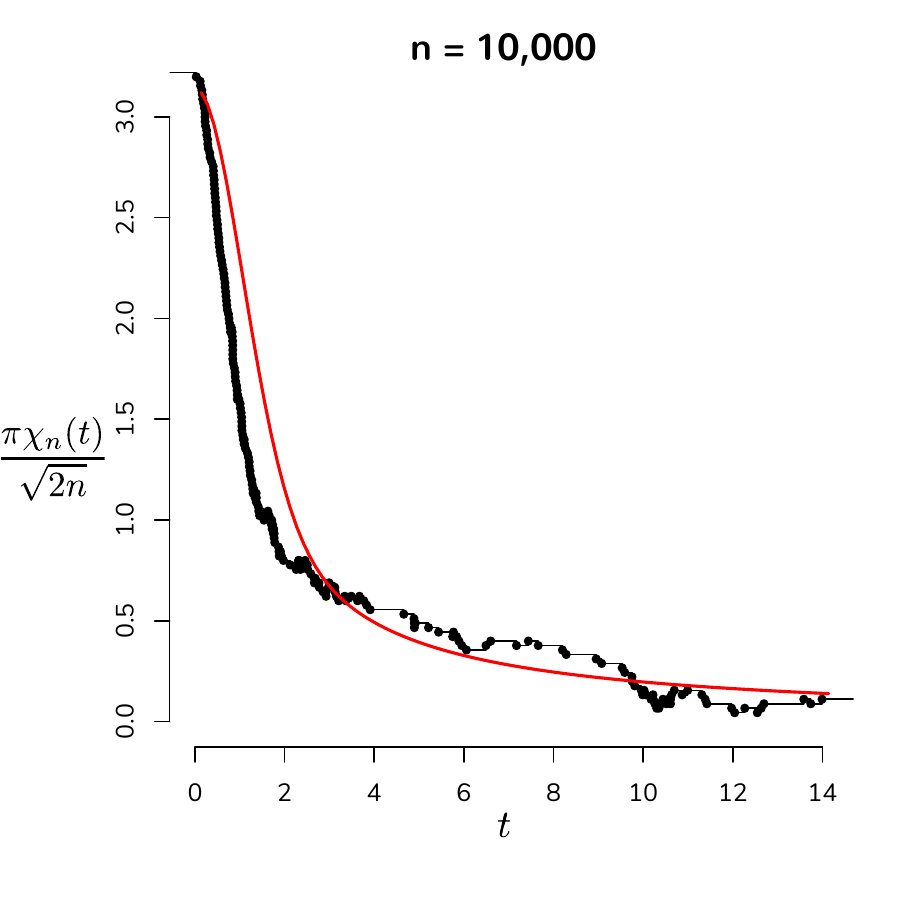}
\includegraphics[width=1.5in]{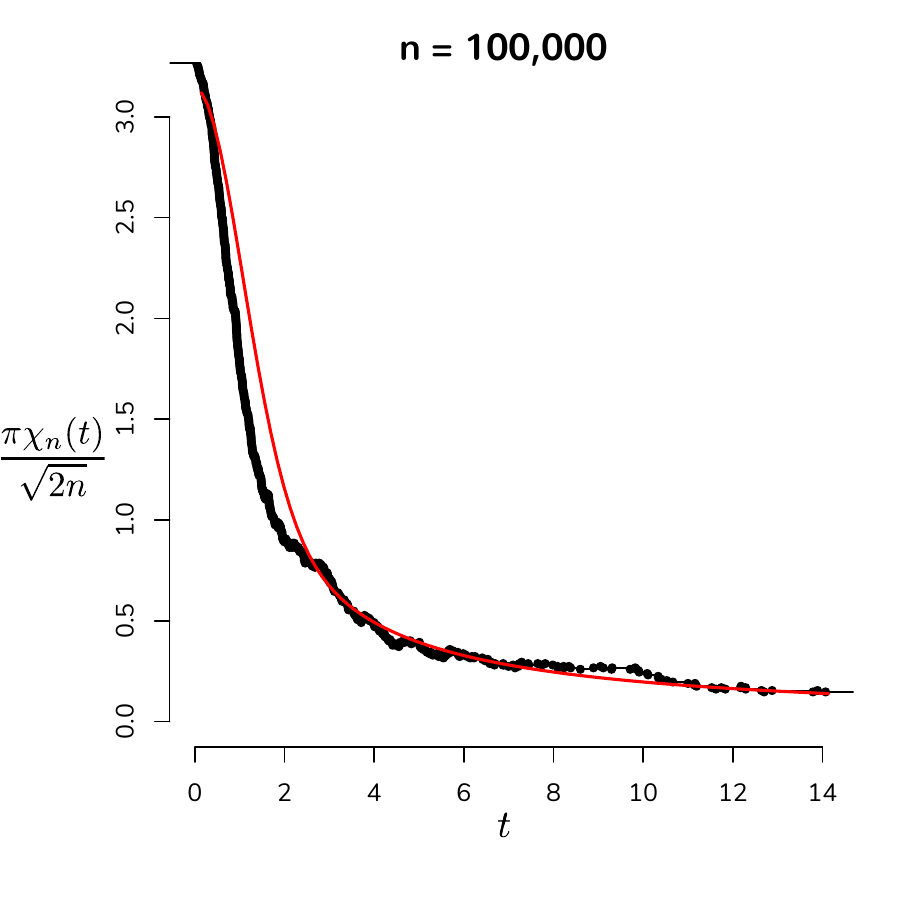}
\caption{\footnotesize{Plots of random realizations of $\sqrt{\pi\omega_d/(2n)}\, \chi_n(t)$ for $d=2$ (in black) in the setup of Example~\ref{ex:chi_np}. In the plots above, as $n$ increases from left to right, the random function converges uniformly to $\sum_{k=0}^{\infty} (-1)^{k}s_{k}(t)$ (in red).}}
\label{fig:chi_np}
\end{figure}
\begin{figure}[t]
\centering
\label{f:evt_ec_proc}
\includegraphics[width = 2.5in]{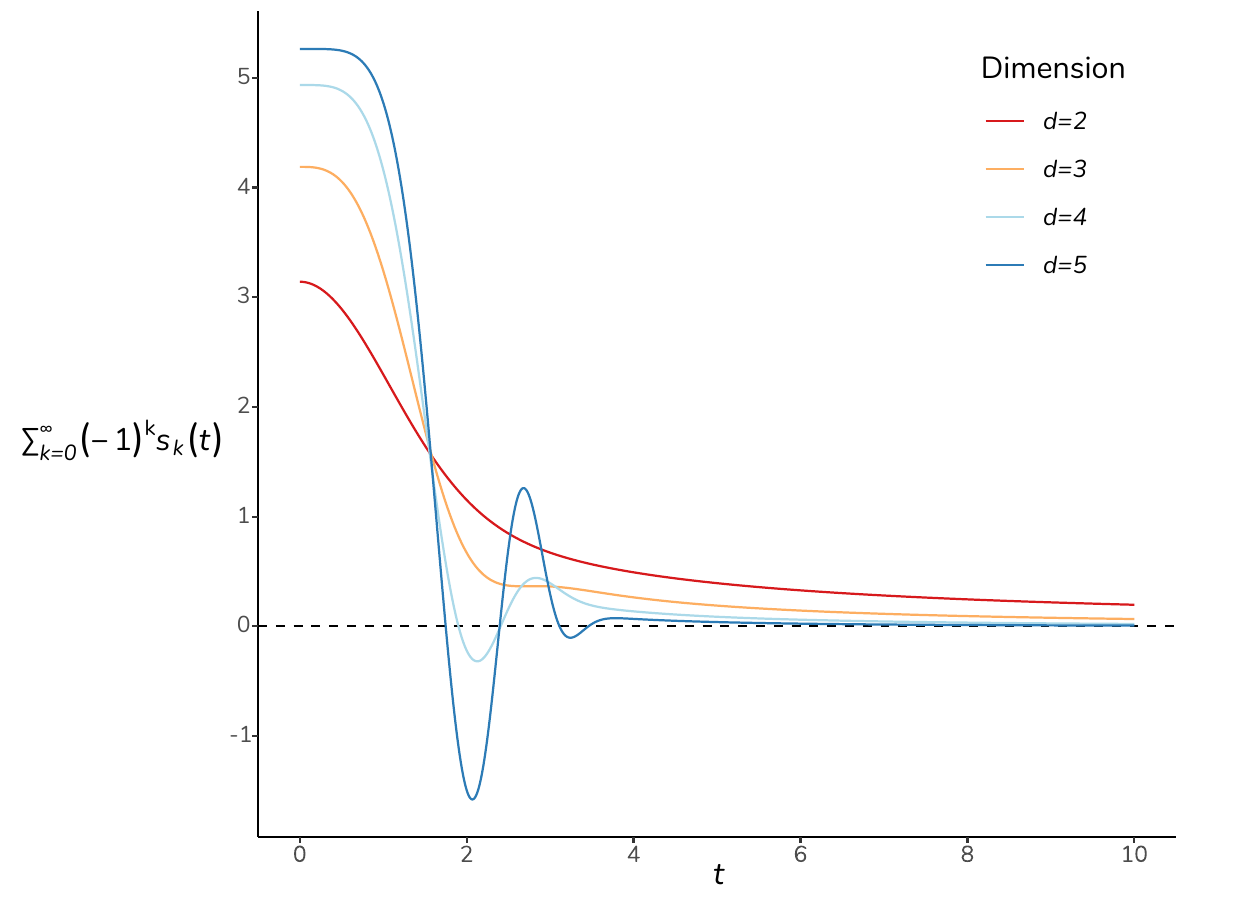}
\caption{Plots of $\sum_{k=0}^{\infty} (-1)^k s_k(t)$ at \eqref{e:ex_dens} for $d=2, 3,4,5$.}
\label{fig:multi.dim.EC}
\end{figure}
 
One of the implications of Theorem \ref{t:heavy} is that one can immediately obtain various limit theorems of functions of the Euler characteristic process. For every continuous function $T$ on $D[0,\infty)$, it indeed holds that, as $n\to\infty$, 
$$
T\left(\sqrt{\frac{\pi\omega_d}{2n}} \chi_n \right) \to T\left(  \sum_{k=0}^\infty (-1)^k s_k \right), \ \ \mathrm{a.s.}
$$
For example, if $U_{a,b}:D[0,\infty) \to [0,\infty)$ is defined by $U_{a,b}(x) := \sup_{a\le x \le b}\big| x(t) \big|$, $0 \le a < b < \infty$, we have
$$
\sqrt{\frac{\pi\omega_d}{2n}} \sup_{a\le t \le b} \big| \chi_n(t) \big| \to \sup_{a\le t \le b} \bigg|  \sum_{k=0}^\infty (-1)^k s_k(t) \bigg| \ \ \mathrm{a.s.}
$$
Furthermore, let $I:D[0, \infty) \to D[0,\infty)$ be defined by $I(x)(t) := \int^{t}_0 x(r) \dif{r}$; then, 
$$
\bigg( \sqrt{\frac{\pi\omega_d}{2n}} \int^t_0 \chi_n(r) \dif{r}, \ t \geq 0\bigg) \to \bigg( \sum_{k=0}^\infty (-1)^k \int^t_0 s_k(r) \dif{r}, \ t \geq 0 \bigg), \ \ \mathrm{a.s. \ in \ } D[0,\infty). 
$$
This result is especially important for applications in TDA. Indeed, $\int^t_0 \chi_n(r) \dif{r}$ represents an alternating sum of the total length of  persistence barcodes of all dimensions, up to time $t$. A persistence barcode is a graphical descriptor of persistent homology, which allows us to visualize the birth time and death time of cycles \citep{ghrist:2008, carlsson:2009}. In light of the TDA literature \citep[e.g., Section 6 of][]{bobrowski2012}, the limit $\int^\infty_0 \chi_n(r) \dif{r}$ is defined as the Euler characteristic of persistence barcodes of the filtration \eqref{e:geom.comp.EVT}. This gives us an estimate of how long the cycles of any dimension live in our extreme sample cloud. 
\end{example}


\section{Exponentially decaying tail case} \label{s:lite}

In this section, we consider a density of an exponentially decaying tail. We assume that the density $f$ is specified by 
\begin{equation} \label{e:edecay}
f(x) = C\exp\big\{-\psi(\,\norm{x})\big\}, \ x\in \R^d,
\end{equation}
where $C$  is a normalizing constant and $\psi: [0,\infty)\to[0,\infty)$ is a regularly varying function (at infinity) of an exponent $\tau \in (0,1]$. Moreover, $\psi$ is assumed to be twice differentiable, such that $\psi^\prime (x)>0$ for all $x>0$, and $\psi^\prime$ is eventually non-increasing. Namely, there exists $z_0>0$ such that $\psi'$ is non-increasing in $[z_0, \infty)$. 
Under this setup, let $a(z):= 1/\psi'(z)$; then, $a$ is also regularly varying with index $1-\tau$ \citep[see, e.g., Proposition 2.5 in][]{htresnick}.  

Here, it is important to note that the occurrence of topological crackle depends on the limit value of $a(z)$ as $z\to\infty$ \citep[see][]{owada_crackle}. In particular, \cite{owada_crackle} showed that crackle occurs if and only if 
\begin{equation} \label{e:crackle.cond}
\zeta:= \lim_{z\to\infty}a(z) \in (0,\infty].
\end{equation}
Since the main theme of this study is topological crackle, we do not treat the case $\zeta=0$. In terms of the regular variation exponent of $\psi$, we exclude the case $\tau>1$. So, for instance, the multivariate Gaussian densities do not belong to the scope of our study. Note that  \eqref{e:crackle.cond} trivially holds for every $\tau \in (0,1)$. 

Now, we are ready to state the FSLLN below. Interestingly, if $\zeta=\infty$ in \eqref{e:crackle.cond}, the limiting function in \eqref{e:sekt} agrees with \eqref{e:limit.heavy} up to multiplicative constants.

\begin{theorem} \label{t:lite}
Suppose that $f$ is a density specified by \eqref{e:edecay} with $\tau\in (0,1]$. If $\tau=1$, we assume \eqref{e:crackle.cond}. If $d=2$, we restrict the range of $\tau$ to $(0,1)$.  
Suppose further that $nf(R_n) \to \xi$ as $n\to\infty$ for some $\xi\in (0,\infty)$. Then, we have, as $n\to\infty$, 
\begin{equation}  \label{e:fslln.exp}
\left( \frac{\chi_n(t)}{a(R_n)R_n^{d-1}}, \, t\ge 0\right) \to \Big( \sum_{k=0}^{\infty} (-1)^{k}s_{k}(t), \, t \ge 0\Big),  \quad \mathrm{a.s. \ in \ } D[0, \infty), 
\end{equation}
where 
\begin{align} 
s_k(t) &:= \frac{\xi^{k+1}}{(k+1)!} \int_0^\infty \int_{S^{d-1}} \int_{(\R^d)^k} \, h_t(0,y_1,\dots,y_k)\, e^{-(k+1)\rho -\zeta^{-1} \sum_{i=1}^k \langle \theta, y_i \rangle}  \label{e:sekt} \\
&\qquad \qquad \qquad \qquad \times \prod_{i=1}^k \one \big\{ \rho+\zeta^{-1} \langle \theta, y_i \rangle \ge 0 \big\}  \dif \y  \dif \nu_{d-1}(\theta) \dif \rho,   \ \ \ t\ge 0, \ \ k \ge 1, \notag
\end{align}
with $\langle \cdot, \cdot \rangle$ being the Euclidean inner product and $s_0(t)\equiv s_{d-1} \xi$. In particular, the limit in \eqref{e:fslln.exp} is convergent for all $t\ge0$. 
\end{theorem}

\begin{example}  \label{ex:exp.density}
We consider a special case of the density in \eqref{e:edecay}, 
$$
f(x) = Ce^{-\|x \|^\tau/\tau}, \ \ \ x\in \R^d, \ \ \tau\in (0,1]. 
$$
Define $R_n = \big( \tau \log n + \tau \log C \big)^{1/\tau}$ so that $nf(R_n) = 1$. Then, $a(z) =z^{1-\tau}$, $z>0$. According to Theorem \ref{t:lite}, 
$$
\left( \frac{\chi_n(t)}{( \tau\log n )^{\frac{d-\tau}{\tau}}}, \, t\ge 0\right) \to \Big( \sum_{k=0}^{\infty} (-1)^{k} s_{k}(t), \, t \ge 0\Big), \ \ \ \mathrm{a.s.\  in\ } D[0,\infty). 
$$
where $s_k(t)$ is defined in \eqref{e:sekt}. Moreover, applying the continuous functions $U_{a,b}$ and $I$ from Example~\ref{ex:chi_np}, we have
 $$
\frac{\sup_{a\le t \le b} \big| \chi_n(t)\big|}{ ( \tau\log n )^{\frac{d-\tau}{\tau}}} \to \sup_{a\le t \le b} \bigg| \sum_{k=0}^{\infty} (-1)^{k}s_{k}(t) \bigg|   \ \ \ \mathrm{a.s.}
$$
and
$$
\Bigg( \frac{\int^t_0 \chi_n(r) \dif{r}}{( \tau\log n )^{\frac{d-\tau}{\tau}}}, \ t \geq 0\Bigg) \to \bigg( \sum_{k=0}^\infty (-1)^k \int^t_0 s_k(r) \dif{r}, \ t \geq 0 \bigg), \ \ \mathrm{a.s. \ in \ } D[0,\infty). 
$$
\end{example}

\section{Proofs of main theorems}  \label{s:proof}

Throughout the proof, denote by $C^*$ a positive constant that is independent of $n$ and may vary between (and even within) the lines. 
Denote by $\text{RV}_\rho$ the collection of regularly varying sequences (or functions) at infinity with exponent $\rho\in\R$. For $a, b \in \R$, write $a\wedge b =\min\{a,b\}$ and $a\vee b=\max\{a,b\}$. For two sequences $(a_n)_{n\geq1}$ and $(b_n)_{n \geq 1}$, $a_n\sim b_n$ means $a_n/b_n\to 1$ as $n\to\infty$. 

First, we present a fundamental result which allows us to extend a pointwise SLLN to a functional SLLN in the space $D[0,\infty)$. 

\begin{proposition}[Proposition 4.2 in \citealp{thomas:owada:2020}] \label{p:monoFSLLN} 
Let $(X_n, \, n \in \mathbb{N})$ be a sequence of random elements in $D[0, \infty)$ with non-decreasing sample paths. Suppose $a: [0,\infty) \to \R$ is a deterministic, continuous, and non-decreasing function. If we have 
$$
X_n(t) \to a(t), \ \    n\to\infty, \ \ \mathrm{a.s.}, 
$$
for every $t\ge0$, then it follows that 
\begin{equation*} 
\sup_{t \in [0,T]} |X_n(t) - a(t)| \to 0, \ \ n\to\infty, \  \ \mathrm{a.s.},  
\end{equation*}
for every $0\le T<\infty$. Hence, it holds that $X_n(t) \to a(t) \ \mathrm{a.s.}$~in $D[0, \infty)$ under the uniform topology.  
\end{proposition}

By virtue of this proposition, for the proof of Theorem \ref{t:heavy} it suffices to show that as $n\to\infty$, 
$$
\frac{\chi_n(t)}{R_n^d} \to \sum_{k=0}^\infty (-1)^k s_k(t), \ \  \text{a.s.}
$$
for every $t\ge 0$. Subsequently, we divide the Euler characteristic process into two terms: 
\begin{equation}  \label{e:decomp.EC.proc}
\chi_n(t) = \sum_{k=0}^\infty S_{2k,n}(t) - \sum_{k=0}^\infty S_{2k+1,n}(t) =: \chi_n^{(1)}(t) - \chi_n^{(2)}(t). 
\end{equation}
In addition, the limiting function can also be decomposed as 
\begin{equation}  \label{e:decomp.limit}
\sum_{k=0}^\infty (-1)^k s_k(t) = \sum_{k=0}^\infty s_{2k}(t)  - \sum_{k=0}^\infty s_{2k+1}(t) =: K_1(t) - K_2(t). 
\end{equation}
From \eqref{e:decomp.EC.proc} and \eqref{e:decomp.limit}, it is now sufficient to prove that for every $t\ge0$ and $i=1,2$, 
\begin{equation}  \label{e:SLLN.EC.decomp.heavy}
\frac{\chi_n^{(i)}(t)}{R_n^d} \to K_i(t), \ \ \ n\to\infty, \ \text{a.s.}
\end{equation}

In the case of Theorem \ref{t:lite}, defining $\chi_n^{(i)}(t)$ and $K_i(t)$  analogously, it suffices to show that for each $t\ge0$ and $i=1,2$, 
\begin{equation}  \label{e:SLLN.EC.decomp.light}
\frac{\chi_n^{(i)}(t)}{a(R_n) R_n^{d-1}} \to K_i(t), \ \ \ n\to\infty, \ \text{a.s.}
\end{equation}


\subsection{Proof of Theorem \ref{t:heavy}}   \label{s:heavy_proofs}

The goal of this subsection is to prove \eqref{e:SLLN.EC.decomp.heavy}. We handle the case $i=1$ only, as the proof is totally the same regardless of $i\in\{ 1,2\}$. Let
\begin{equation}  \label{e:def.vm}
u_m = \lfloor e^{m^\gamma} \rfloor, \ \ \ m=0,1,2,\dots,
\end{equation}
for some $\gamma \in (0,1)$. Then, for every $n\in \N$, there exists a unique $m=m(n)$ such that $u_m \le n < u_{m+1}$. Let us also define 
\begin{align}
p_m &= \text{argmax} \{ u_m \le \ell \le u_{m+1}: R_\ell \}, \label{e:def.pm} \\
q_m &= \text{argmin} \{ u_m \le \ell \le u_{m+1}: R_\ell \}. \label{e:def.qm}
\end{align}
It then holds that $R_{p_m} = \max_{v_m \le \ell \le v_{m+1}}R_\ell$ and $R_{q_m} = \min_{v_m \le \ell \le v_{m+1}}R_\ell$.

Below, we offer a lemma on the asymptotic moments of certain variants of the process $\chi_n^{(1)}(t)$, defined by 
\begin{align}
T_m(t) &:= \sum_{k=0}^\infty \sum_{\substack{\ \Y\subset \X_{u_{m+1}}, \\ |\Y|=2k+1}} h_t(\Y)\, \one \big\{  \min_{y\in \Y}\|y\|\ge R_{q_m} \big\}, \label{e:def.Tm}\\
U_m(t) &:= \sum_{k=0}^\infty \sum_{\substack{\ \Y\subset \X_{u_{m}}, \\ |\Y|=2k+1}} h_t(\Y)\, \one \big\{  \min_{y\in \Y}\|y\|\ge R_{p_m} \big\}. \label{e:def.Um}
\end{align} 

\begin{lemma}  \label{l:var_res}
Under the assumptions of Theorem \ref{t:heavy}, we have the following asymptotic results on the first and second moments of $T_m(t)$ and $U_m(t)$. 
\begin{align}
&\lim_{m\to\infty} R_{q_m}^{-d} \E \big[ T_m (t)\big] = K_1(t), \label{e:first.Tm} \\
&\lim_{m\to\infty} R_{p_m}^{-d} \E \big[ U_m (t)\big] = K_1(t), \label{e:first.Um} \\
&\sup_{m\ge1} R_{q_m}^{-d} \Var\big( T_m(t) \big) < \infty, \label{e:second.Tm} \\
&\sup_{m\ge 1} R_{p_m}^{-d} \Var\big( U_m(t) \big) < \infty. \label{e:second.Um}
\end{align}
\end{lemma}

\begin{proof}
We begin by offering the proofs of \eqref{e:first.Tm} and \eqref{e:first.Um} by extending the argument in the proof of Proposition 7.2 of \cite{owada2017}. As for \eqref{e:first.Tm}, it is clear that 
\begin{align*}
R_{q_m}^{-d} \E \big[ T_m (t)\big] = \sum_{k=0}^\infty R_{q_m}^{-d} \binom{u_{m+1}}{2k+1} \E\Big[  h_t(X_1,\dots,X_{2k+1})\, \one \big\{\min_{1\le i \le 2k+1} \|X_i\| \ge R_{q_m} \big\} \Big], 
\end{align*}
where $X_1,\dots,X_{2k+1}$ are i.i.d random variables with density $f$. From this, we have 
\begin{align}
&R_{q_m}^{-d} \binom{u_{m+1}}{2k+1} \E\Big[  h_t(X_1,\dots,X_{2k+1})\, \one \big\{ \min_{1\le i \le 2k+1}\|X_i\| \ge R_{q_m} \big\} \Big]  \label{e:first.change.variables}\\
&= R_{q_m}^{-d} \binom{u_{m+1}}{2k+1} \int_{(\R^d)^{2k+1}} h_t(x_1,\dots,x_{2k+1}) \prod_{i=1}^{2k+1} f(x_i)\one  \big\{ \|x_i\| \ge R_{q_m} \big\} \dif \x \notag \\
&= R_{q_m}^{-d} \binom{u_{m+1}}{2k+1} \int_{\R^d} \int_{(\R^d)^{2k}} h_t(0,y_1,\dots,y_{2k}) f(x) \one \big\{ \|x\| \ge R_{q_m} \big\} \notag \\
&\qquad \qquad \qquad \qquad \qquad \qquad \qquad  \times \prod_{i=1}^{2k} f(x+y_i) \one \big\{ \|x+y_i\| \ge R_{q_m}\big\} \dif \y \dif x, \notag
\end{align}
by the change of variables $x_i=x+y_{i-1}$, $i=1,\dots,2k+1$ (with $y_0\equiv 0$) and the translation invariance of $h_t$. Furthermore, we make the change of variables by $x=R_{q_m}\rho\theta$ with $\rho\ge 1$ and $\theta\in S^{d-1}$, to get that 
\begin{align}
&R_{q_m}^{-d} \binom{u_{m+1}}{2k+1} \E\Big[  h_t(X_1,\dots,X_{2k+1})\, \one \big\{ \min_{1\le i \le 2k+1}\|X_i\| \ge R_{q_m} \big\} \Big] \label{e:after.polar} \\
&= \binom{u_{m+1}}{2k+1} f(R_{q_m})^{2k+1} \int_1^\infty \int_{S^{d-1}} \int_{(\R^d)^{2k}} h_t(0,\y) \rho^{d-1} \frac{f(R_{q_m}\rho)}{f(R_{q_m})} \notag \\
&\qquad \qquad \qquad \qquad \qquad  \times \prod_{i=1}^{2k} \frac{f\big( R_{q_m} \|\rho\theta +y_i/R_{q_m}\| \big)}{f(R_{q_m})}\, \one \big\{ \|\rho\theta + y_i /R_{q_m}\| \ge 1 \big\} \dif \y \dif \nu_{d-1}(\theta) \dif \rho, \notag
\end{align}
where $\y=(y_1,\dots,y_{2k})\in (\R^d)^{2k}$. 
Next, for a fixed constant $\eta\in (0,\alpha-d)$, Potter's bounds \citep[see Proposition 2.6 in][]{htresnick} yield that 
\begin{equation}  \label{e:Potter1}
\frac{f(R_{q_m}\rho)}{f(R_{q_m})} \le 2 \rho^{-\alpha +\eta}, 
\end{equation}
and 
\begin{equation}  \label{e:Potter2}
\prod_{i=1}^{2k} \frac{f\big( R_{q_m} \|\rho\theta +y_i/R_{q_m}\| \big)}{f(R_{q_m})}\, \one \big\{ \|\rho\theta + y_i /R_{q_m}\| \ge 1 \big\} \le 2^{2k}
\end{equation}
for sufficiently large $m$. Since $\int_1^\infty \rho^{d-1-\alpha +\eta}\dif \rho<\infty$ and $\int_{(\R^d)^{2k}}h_t(0,\y)\dif \y<\infty$ by property \textbf{(H3)} of $h$, we can see that the regular variation of $f$, as well as the dominated convergence theorem, ensures that the triple integral in \eqref{e:after.polar} converges to 
$$
s_{d-1} \int_1^\infty \rho^{d-1-\alpha (2k+1)} \dif \rho \int_{(\R^d)^{2k}} h_t(0,\y) \dif \y = \frac{s_{d-1}}{\alpha(2k+1)-d}\, \int_{(\R^d)^{2k}} h_t(0,\y) \dif \y. 
$$
Furthermore, \eqref{e:weak.core.regime} ensures that as $m\to\infty$, 
$$
\binom{u_{m+1}}{2k+1} f(R_{q_m})^{2k+1} \sim \frac{\big( u_{m+1}f(R_{q_m})\big)^{2k+1}}{(2k+1)!} \sim \frac{\xi^{2k+1}}{(2k+1)!}\, \Big(  \frac{u_{m+1}}{q_m}\Big)^{2k+1}, 
$$
so that 
\begin{equation}  \label{e:ratio.um.qm}
1 \le \frac{u_{m+1}}{q_m} \le \frac{u_{m+1}}{u_m} \le \frac{e^{(m+1)^\gamma - m^\gamma}}{1-e^{-m^{\gamma}}} = \frac{e^{m^{\gamma-1}(\gamma + o(1))}}{1-e^{-m^{\gamma}}}. 
\end{equation}
As $0<\gamma <1$, the rightmost term above converges to $1$ as $m\to\infty$. Hence, 
$$
\binom{u_{m+1}}{2k+1} f(R_{q_m})^{2k+1} \to \frac{\xi^{2k+1}}{(2k+1)!}, \ \ \text{as } m\to\infty. 
$$
Combining all of these results, it follows that 
$$
R_{q_m}^{-d} \binom{u_{m+1}}{2k+1} \E\Big[  h_t(X_1,\dots,X_{2k+1})\, \one \big\{ \min_{1\le i \le 2k+1}\|X_i\| \ge R_{q_m} \big\} \Big]  \to s_{2k}(t), \ \ \text{as } m\to\infty, 
$$
which yields \eqref{e:first.Tm} as desired. 
The proof of \eqref{e:first.Um} is almost the same, so we omit it here. 

Now we will prove \eqref{e:second.Tm}. We see that 
\begin{align*}
\E \big[ T_m(t)^2 \big] &=\sum_{k_1=0}^\infty \sum_{k_2=0}^\infty \sum_{\ell=0}^{2 (k_1 \wedge k_2)+1} \E \bigg[ \prod_{i=1}^2 \bigg( \sum_{\substack{\Y_i\subset \X_{u_{m+1}}, \\  |\Y_i|=2k_i+1}} h_t(\Y_i)\,  \one \big\{  \min_{y\in\Y_i} \|y\| \ge R_{q_m} \big\} \bigg)\, \one \big\{  |\Y_1 \cap \Y_2| =\ell \big\}\bigg] \\
&= \sum_{k_1=0}^\infty \sum_{k_2=0}^\infty \sum_{\ell=0}^{2 (k_1 \wedge k_2)+1} \binom{u_{m+1}}{2(k_1+k_2)+2-\ell} \frac{\big( 2(k_1+k_2)+2-\ell \big)!}{(2k_1+1-\ell)!\, (2k_2+1-\ell)!\, \ell!}\,  \\
&\times \E \Big[  h_t(X_1,\dots,X_{2k_1+1}) h_t (X_1,\dots,X_\ell, X_{2k_1+2}, \dots, X_{2(k_1+k_2)+2-\ell}) \\
&\qquad \qquad \qquad\qquad \qquad \qquad\qquad \qquad \qquad\times  \one \big\{ \min_{1\le i \le 2(k_1+k_2)+2-\ell}\|X_i\|\ge R_{q_m} \big\} \Big] \\
&=: \sum_{k_1=0}^\infty \sum_{k_2=0}^\infty \sum_{\ell=0}^{2 (k_1 \wedge k_2)+1} \binom{u_{m+1}}{2(k_1+k_2)+2-\ell} \frac{\big( 2(k_1+k_2)+2-\ell \big)!}{(2k_1+1-\ell)!\, (2k_2+1-\ell)!\, \ell!}\, \E[I_{k_1,k_2,\ell}], 
\end{align*}
where $X_1,\dots,X_{2(k_1+k_2)+2-\ell}$ are i.i.d random points with density $f$. In the above, if $\ell=0$, we take that 
$$
h_t (X_1,\dots,X_\ell, X_{2k_1+2}, \dots, X_{2(k_1+k_2)+2-\ell}) := h_t(X_{2k_1+2}, \dots, X_{2(k_1+k_2)+2}). 
$$
From this, we see that 
\begin{align}
&R_{q_m}^{-d} \text{Var} \big( T_m(t) \big) \label{e:scaled.variance.Tm}\\
&= \sum_{k_1=0}^\infty \sum_{k_2=0}^\infty \sum_{\ell=1}^{2 (k_1 \wedge k_2)+1} R_{q_m}^{-d}  \binom{u_{m+1}}{2(k_1+k_2)+2-\ell} \frac{\big( 2(k_1+k_2)+2-\ell \big)!}{(2k_1+1-\ell)!\, (2k_2+1-\ell)!\, \ell!}\, \E[I_{k_1,k_2,\ell}] \notag \\
&\quad + \sum_{k_1=0}^\infty \sum_{k_2=0}^\infty R_{q_m}^{-d}  \binom{u_{m+1}}{2(k_1+k_2)+2} \binom{2(k_1+k_2)+2}{2k_1+1} \E [I_{k_1,k_2,0}] \notag \\
&\quad -  \sum_{k_1=0}^\infty \sum_{k_2=0}^\infty R_{q_m}^{-d}  \binom{u_{m+1}}{2k_1+1} \binom{u_{m+1}}{2k_2+1} \E [I_{k_1,k_2,0}] \notag \\
&\le \sum_{k_1=0}^\infty \sum_{k_2=0}^\infty \sum_{\ell=1}^{2 (k_1 \wedge k_2)+1} R_{q_m}^{-d}  \binom{u_{m+1}}{2(k_1+k_2)+2-\ell} \frac{\big( 2(k_1+k_2)+2-\ell \big)!}{(2k_1+1-\ell)!\, (2k_2+1-\ell)!\, \ell!}\, \E[I_{k_1,k_2,\ell}], \notag
\end{align}
where the last inequality comes from 
\begin{align*}
&\binom{u_{m+1}}{2(k_1+k_2)+2} \binom{2(k_1+k_2)+2}{2k_1+1}  - \binom{u_{m+1}}{2k_1+1} \binom{u_{m+1}}{2k_2+1} \\
&\qquad = \dbinom{u_{m+1}}{2k_1+1}\Bigg[ \dbinom{u_{m+1} - 2k_1-1}{2k_2+1} - \dbinom{u_{m+1}}{2k_2+1}\Bigg] < 0.
\end{align*}

For our purposes, we must examine the appropriate upper bounds of $\E[I_{k_1,k_2,\ell}]$ for $\ell\ge1$. For every $\ell\ge 1$, performing the change of variables, $x_i=x+y_{i-1}$, $i=1,\dots,2(k_1+k_2)+2-\ell$ (with $y_0\equiv 0$), we have that 
\begin{align*}
\E[I_{k_1,k_2,\ell}] &= \int_{(\R^d)^{2(k_1+k_2)+2-\ell}} \hspace{-10pt} h_t(x_1,\dots,x_{2k_1+1})h_t(x_1,\dots,x_\ell, x_{2k_1+2}, \dots, x_{2(k_1+k_2)+2-\ell}) \\
&\qquad \qquad \qquad \qquad \qquad \qquad \qquad\times \prod_{i=1}^{2(k_1+k_2)+2-\ell} f(x_i) \one \big\{ \|x_i\|\ge R_{q_m} \big\} \dif \x \\
&= \int_{\R^d}\int_{(\R^d)^{2(k_1+k_2)+1-\ell}} \hspace{-10pt} h_t(0,y_1,\dots,y_{2k_1})h_t(0,y_1,\dots,y_{\ell-1}, y_{2k_1+1}, \dots, y_{2(k_1+k_2)+1-\ell}) \\
&\qquad \qquad \qquad\times f(x)\, \one \big\{ \|x\|\ge R_{q_m} \big\} \prod_{i=1}^{2(k_1+k_2)+1-\ell} f(x+y_i)\, \one \big\{ \|x+y_i\| \ge R_{q_m} \big\} \dif \y \dif x. 
\end{align*}
As in \eqref{e:after.polar}, we apply the polar coordinate transform $x=R_{q_m}\rho\theta$ with $\rho\ge 1$ and $\theta\in S^{d-1}$, to obtain that 
\begin{align*}
&\E[I_{k_1,k_2,\ell}] = R_{q_m}^d f(R_{q_m})^{2(k_1+k_2)+2-\ell} \int_1^\infty \int_{S^{d-1}} \int_{(\R^d)^{2(k_1+k_2)+1-\ell}} \hspace{-10pt} h_t(0,\y_0, \y_1)h_t(0,\y_0,\y_2) \rho^{d-1}\\
& \times \frac{f(R_{q_m}\rho)}{f(R_{q_m})} \prod_{i=1}^{2(k_1+k_2)+1-\ell} \frac{f\big( R_{q_m} \| \rho\theta +y_i /R_{q_m} \| \big)}{f(R_{q_m})} \, \ind{\| \rho\theta +y_i /R_{q_m} \| \geq 1} \dif\, (\y_0\cup \y_1 \cup \y_2)\dif \nu_{d-1}(\theta) \dif \rho, 
\end{align*}
where $\y_0 = (y_1, \dots, y_{\ell-1})$, $\y_1 = (y_\ell, \dots, y_{2k_1})$ and $\y_2 = (y_{2k_1+1}, \dots, y_{2(k_1+k_2)+1-\ell})$. 
Appealing to Potter's bounds as in \eqref{e:Potter1} and \eqref{e:Potter2}, as well as \eqref{e:weak.core.regime}, there exists $N\in\N$ such that for all $m\ge N$, 
\begin{align}
&R_{q_m}^{-d} \binom{u_{m+1}}{2(k_1+k_2)+2-\ell} \E[I_{k_1,k_2,\ell}] \notag \\
&\le \frac{\big(u_{m+1}f(R_{q_m})  \big)^{2(k_1+k_2)+2-\ell}}{\big( 2(k_1+k_2)+2-\ell \big)!}\, \int_1^\infty \int_{S^{d-1}} \int_{(\R^d)^{2(k_1+k_2)+1-\ell}} \hspace{-10pt} h_t(0,\y_0, \y_1)h_t(0,\y_0,\y_2) \rho^{d-1} \notag \\
&\qquad \qquad \qquad \qquad \qquad \qquad \times 2\rho^{-\alpha +\eta} \times 2^{2(k_1+k_2)+1-\ell} \dif\,  (\y_0\cup \y_1 \cup \y_2) \dif \nu_{d-1}(\theta) \dif \rho \notag \\
&\le \frac{(4 \xi)^{2(k_1+k_2)+2-\ell}}{\big( 2(k_1+k_2)+2-\ell \big)!} \cdot \frac{s_{d-1}}{\alpha -d -\eta} \int_{(\R^d)^{2(k_1+k_2)+1-\ell}} \hspace{-10pt} h_t(0,\y_0, \y_1)h_t(0,\y_0,\y_2) \dif\,  (\y_0\cup \y_1 \cup \y_2).  \notag
\end{align}
By virtue of property \textbf{(H3)} of $h$, 
\begin{align*}
\int_{(\R^d)^{2(k_1+k_2)+1-\ell}} \hspace{-10pt} h_t(0,\y_0, \y_1)h_t(0,\y_0,\y_2) \dif\,  (\y_0\cup \y_1 \cup \y_2) &\le \int_{(\R^d)^{2(k_1+k_2)+1-\ell}} \hspace{-10pt} \prod_{i=1}^{2(k_1+k_2)+1-\ell} \hspace{-10pt} \one \big\{ \|y_i\| \le ct \big\} \dif \y \\
&= \big((ct)^d\omega_d\big)^{2(k_1+k_2)+1-\ell}. 
\end{align*}
Therefore, for all $m\ge N$, 
\begin{align}
R_{q_m}^{-d} \binom{u_{m+1}}{2(k_1+k_2)+2-\ell} \E[I_{k_1,k_2,\ell}] \notag &\le \frac{(4 \xi)^{2(k_1+k_2)+2-\ell} \big((ct)^d\omega_d\big)^{2(k_1+k_2)+1-\ell} }{\big( 2(k_1+k_2)+2-\ell \big)!} \cdot \frac{s_{d-1}}{\alpha -d -\eta} \notag \\
&\le \frac{(C^*)^{2(k_1+k_2)+2-\ell} }{\big( 2(k_1+k_2)+2-\ell \big)!}. \label{e:bound.EIkl}
\end{align}
Note that the constant $C^*$ does not depend on $k_1$, $k_2$, and $\ell$. 
Returning to \eqref{e:scaled.variance.Tm} and applying the bound in \eqref{e:bound.EIkl}, we have that 
\begin{align*}
R_{q_m}^{-d} \text{Var} \big( T_m(t) \big) &\le \sum_{k_1=0}^\infty \sum_{k_2=0}^\infty \sum_{\ell=1}^{2 (k_1 \wedge k_2)+1}  \frac{(C^*)^{2(k_1+k_2)+2-\ell}}{(2k_1+1-\ell)!\,(2k_2+1-\ell)!\, \ell !} \\
&\le 2 \sum_{k_1=0}^\infty \sum_{k_2=0}^{k_1} \sum_{\ell=1}^{k_2+1}  \frac{(C^*)^{k_1+k_2+2-\ell}}{(k_1+1-\ell)!\,(k_2+1-\ell)!\, \ell !} \\
&\le 2 \sum_{\ell=1}^\infty \sum_{k_1=\ell-1}^\infty  \sum_{k_2=\ell-1}^\infty \frac{(C^*)^{k_1+1-\ell}}{(k_1+1-\ell)!} \cdot  \frac{(C^*)^{k_2+1-\ell}}{(k_2+1-\ell)!} \cdot \frac{(C^*)^\ell}{\ell !} = 2e^{3C^*} <\infty. 
\end{align*}
Since the proof of \eqref{e:second.Um} is very similar to that of \eqref{e:second.Tm}, we will omit it.
\end{proof}

\begin{proof}[Theorem \ref{t:heavy}]
By the definition of $u_m$, $p_m$, and $q_m$, we have, for every $n\ge 1$, 
$$
\frac{U_m(t)}{R_{p_m}^d} \le \frac{\chi_n^{(1)}(t)}{R_n^d} \le \frac{T_m(t)}{R_{q_m}^d}. 
$$
Then, Lemma \ref{l:var_res} gives that 
\begin{equation} \label{e:T_m.bound}
\limsup_{n \to \infty}  \frac{\chi_n^{(1)}(t)}{R_n^d} \leq K_1(t) + \limsup_{m \to \infty} \frac{T_m(t)-\E[T_m(t)]}{R_{q_m}^d}, \ \ \text{a.s.}, 
\end{equation}
and 
\begin{equation} \label{e:U_m.bound}
\liminf_{n \to \infty} \frac{\chi_n^{(1)}(t)}{R_n^d} \ge  K_1(t) +  \liminf_{m \to \infty} \frac{U_m(t) - \E[U_m(t)]}{R_{p_m}^d}, \ \ \text{a.s.}
\end{equation}
Let us continue by showing that 
\begin{equation}  \label{e:conv.for.T_m}
R_{q_m}^{-d} \big( T_m(t)-\E[T_m(t)] \big) \to 0, \ \ m\to\infty, \ \text{a.s.}, 
\end{equation}
and  
\begin{equation}  \label{e:conv.for.U_m}
R_{p_m}^{-d} \big( U_m(t)-\E[U_m(t)] \big) \to 0, \ \ m\to\infty, \ \text{a.s.}
\end{equation}
For \eqref{e:conv.for.T_m}, it follows from \eqref{e:second.Tm} in Lemma \ref{l:var_res} and Chebyshev's inequality that, for every $\epsilon>0$, 
\begin{align*}
\sum_{m=1}^\infty \P \Big( \big| \, T_m(t) - \E[T_m(t)] \, \big| > \epsilon R_{q_m}^d \Big) \le \sum_{m=1}^\infty \frac{\text{Var}\big(T_m(t)  \big)}{\epsilon^2 (R_{q_m}^d)^2} \le C^* \sum_{m=1}^\infty \frac{1}{R_{q_m}^d}. 
\end{align*}
As $R_n\in \text{RV}_{1/\alpha}$ (see \eqref{e:explicit.Rn}), we have that 
$$
R_{q_m} \ge C^* q_m^{1/(2\alpha)} \ge C^* u_m^{1/(2\alpha)} \ge C^* e^{m^\gamma/(3\alpha)}
$$  
for all $m\ge 1$. Now, we have $\sum_m R_{q_m}^{-d} \le C^* \sum_m e^{-dm^\gamma/(3\alpha)} < \infty$, and the Borel-Cantelli lemma concludes \eqref{e:conv.for.T_m}. The proof of \eqref{e:conv.for.U_m} is analogous by virtue of \eqref{e:second.Um} in Lemma \ref{l:var_res}. Now, combining \eqref{e:T_m.bound},  \eqref{e:U_m.bound}, \eqref{e:conv.for.T_m}, and \eqref{e:conv.for.U_m} completes the proof. 

Finally, let us explicitly demonstrate that the limit in \eqref{e:fslln.heavy} is finite for all $t\ge0$. By virtue of property \textbf{(H3)} of $h$, 
\begin{align*}
\Big| \,  \sum_{k=0}^\infty (-1)^k s_k(t) \, \Big| &\le \sum_{k=0}^\infty \frac{s_{d-1}\xi^{k+1}}{(k+1)! \big( \alpha (k+1)-d \big)}\int_{(\R^d)^k} \prod_{i=1}^k \one \big\{ \|y_i\| \le ct \big\} \dif \y \\
&\le C^* \sum_{k=0}^\infty \frac{\big( (ct)^d \xi \omega_d \big)^k}{k!} = C^* e^{(ct)^d \xi \omega_d } <\infty. 
\end{align*}
\end{proof}

\subsection{Proof of Theorem \ref{t:lite}} \label{s:lite_proofs}

The goal here is to prove \eqref{e:SLLN.EC.decomp.light}. Once again, we deal with the case $i=1$ only. The proof is essentially similar in character to the proof of Theorem \ref{t:heavy} but involves more complex machinery. 
First we take 
\begin{equation}  \label{e:constraint.gamma}
\gamma \in \Big( \frac{\tau}{d-\tau}, 1 \Big)
\end{equation}
(recall that we have restricted the range of $\tau$ to $(0,1)$ when $d=2$), and define 
$$
u_m = \lfloor e^{m^\gamma} \rfloor,  \ \ m=0,1,2,\dots
$$
as in \eqref{e:def.vm}. Moreover, let $p_m$ and $q_m$ remain as before---see \eqref{e:def.pm} and \eqref{e:def.qm}. 
Additionally, we also introduce 
\begin{align*}
v_m &:= \text{argmax} \big\{ u_m \le \ell \le u_{m+1}: a(R_\ell)R_\ell^{d-1} \big\}, \\
w_m &:= \text{argmin} \big\{ u_m \le \ell \le u_{m+1}: a(R_\ell)R_\ell^{d-1} \big\}. 
\end{align*}
Let $T_m(t)$ and $U_m(t)$ be defined in the same way as \eqref{e:def.Tm} and \eqref{e:def.Um}. 

The lemma below is analogous to Lemma \ref{l:var_res} \citep[see also Proposition 7.4 in][]{owada2017} that provides the asymptotic moments of $T_m(t)$ and $U_m(t)$. 

\begin{lemma}  \label{l:var_res_lite}
Under the assumptions of Theorem \ref{t:lite}, we have the following. 
\begin{align}
&\lim_{m\to\infty} \big[ a(R_{w_m})R_{w_m}^{d-1} \big]^{-1} \E \big[ T_m (t)\big] = K_1(t), \label{e:first.Tm.light} \\
&\lim_{m\to\infty} \big[ a(R_{v_m})R_{v_m}^{d-1} \big]^{-1}  \E \big[ U_m (t)\big] = K_1(t), \label{e:first.Um.light} \\
&\sup_{m\ge1} \big[ a(R_{w_m})R_{w_m}^{d-1} \big]^{-1} \Var\big( T_m(t) \big) < \infty, \label{e:second.Tm.light} \\
&\sup_{m\ge 1} \big[ a(R_{v_m})R_{v_m}^{d-1} \big]^{-1}  \Var\big( U_m(t) \big) < \infty. \label{e:second.Um.light}
\end{align}
\end{lemma}

\begin{proof}
We begin by proving \eqref{e:first.Tm.light} and \eqref{e:first.Um.light}. By the same change of variables as in \eqref{e:first.change.variables} and the translation invariance of $h_t$, 
\begin{align*}
&\big[ a(R_{w_m})R_{w_m}^{d-1} \big]^{-1} \E \big[ T_m (t)\big] \\
&=\sum_{k=0}^\infty \big[ a(R_{w_m})R_{w_m}^{d-1} \big]^{-1} \binom{u_{m+1}}{2k+1} \int_{(\R^d)^{2k+1}} h_t(x_1,\dots,x_{2k+1}) \prod_{i=1}^{2k+1} f(x_i)\, \one \big\{ \|x_i\| \ge R_{q_m} \big\} \dif \x \\
&= \sum_{k=0}^\infty \big[ a(R_{w_m})R_{w_m}^{d-1} \big]^{-1} \binom{u_{m+1}}{2k+1} \int_{\R^d} \int_{(\R^d)^{2k}} h_t(0,\y) f(x)\, \one \big\{ \|x\|\ge R_{q_m} \big\} \\
&\qquad \qquad \qquad \qquad \qquad \times \prod_{i=1}^{2k} f(x+y_i)\, \one \big\{ \|x+y_i\| \ge R_{q_m} \big\} \dif \y \dif x, 
\end{align*}
where $\y=(y_1,\dots,y_{2k})\in (\R^d)^{2k}$. 
Here, we make the change of variable by $x=\big( R_{q_m} +a(R_{q_m})\rho\big)\theta$ with $\rho\ge 0$ and $\theta \in S^{d-1}$. Then, the integral above becomes 
\begin{align*}
&a(R_{q_m}) \int_0^\infty \int_{S^{d-1}} \int_{(\R^d)^{2k}} h_t(0,\y) \big( R_{q_m} +a(R_{q_m})\rho\big)^{d-1} f \big( R_{q_m} +a(R_{q_m})\rho\big)  \\
&\qquad \qquad \times \prod_{i=1}^{2k} f\Big( \big\| \big(\stan\big)\theta +y_i \big\| \Big)\, \one \Big\{ \big\| \big(\stan\big)\theta +y_i \big\| \ge R_{q_m} \Big\} \dif \y \dif \nu_{d-1}(\theta) \dif\rho, 
\end{align*}
which implies that 
\begin{align}
&\big[ a(R_{w_m})R_{w_m}^{d-1} \big]^{-1} \E \big[ T_m (t)\big] \label{e:ETm.normalized.light} \\
&= \sum_{k=0}^\infty  \frac{a(R_{q_m})R_{q_m}^{d-1}}{a(R_{w_m})R_{w_m}^{d-1}}\, \binom{u_{m+1}}{2k+1} f(R_{q_m})^{2k+1} \notag  \\
&\qquad  \times \int_0^\infty \int_{S^{d-1}} \int_{(\R^d)^{2k}} h_t(0,\y) \Big( 1+\frac{a(R_{q_m})}{R_{q_m}} \rho \Big)^{d-1} \frac{f\big( \stan \big)}{f(R_{q_m})} \notag \\
&\qquad \quad  \times \prod_{i=1}^{2k} \frac{f\Big( \big\| \big(\stan\big)\theta +y_i \big\| \Big)}{f(R_{q_m})}\, \one \Big\{ \big\| \big(\stan\big)\theta +y_i \big\| \ge R_{q_m} \Big\} \dif \y \dif \nu_{d-1}(\theta) \dif\rho. \notag
\end{align}
For the last expression, we claim that 
\begin{align}
&\frac{a(R_{q_m})R_{q_m}^{d-1}}{a(R_{w_m})R_{w_m}^{d-1}} \to 1, \ \ \ m\to\infty,  \label{e.ratio.aR.R}
\end{align}
and
\begin{align}
&\binom{u_{m+1}}{2k+1} f(R_{q_m})^{2k+1} \to \frac{\xi^{2k+1}}{(2k+1)!},  \ \ \ m\to\infty. \label{e:weak.core.cond.variant}
\end{align}
Because of \eqref{e:weak.core.regime}, we have $\psi(R_n) \sim \log(Cn/\xi)$ as $n\to\infty$. With the assumptions on the density \eqref{e:edecay}, Proposition 2.6 in \cite{htresnick} gives that $\psi^\leftarrow \in \text{RV}_{1/\tau}$. It follows from the uniform convergence of regularly varying sequences \citep[see Proposition 2.4 in ][]{htresnick} that 
$$
\frac{\psi^\leftarrow \big( \psi(R_n) \big)}{\psi^\leftarrow \big( \log (Cn/\xi) \big)} \sim \Big( \frac{\psi(R_n)}{\log (Cn/\xi)} \Big)^{1/\tau} \to 1,  \ \ \text{as } n\to\infty. 
$$
Since $\psi^\leftarrow\big(\psi(R_n)  \big) \sim R_n$ as $n\to\infty$, the above relation and $\log (Cn/\xi) \sim \log n$, $n\to\infty$, implies 
\begin{equation}  \label{e:R_n.asym}
R_n \sim \psi^\leftarrow\big( \log n \big), \ \  \ n\to\infty. 
\end{equation} 

Now we are ready to prove \eqref{e.ratio.aR.R}. By the uniform convergence of regularly varying sequences, 
\begin{equation}  \label{e:ratio.Rqm.Rwm}
\frac{R_{q_m}}{R_{w_m}} \sim \frac{\psi^\leftarrow \big( \log q_m \big)}{\psi^\leftarrow \big( \log w_m \big)} \sim \bigg( \frac{\log q_m}{\log w_m} \bigg)^{1/\tau}, \  \ \ m\to\infty. 
\end{equation}
Notice that 
$$
\frac{u_m}{u_{m+1}} \le \frac{q_m}{w_m} \le \frac{u_{m+1}}{u_m}, 
$$
so that $u_{m+1}/u_m\to1$ as $m\to\infty$ (see \eqref{e:ratio.um.qm}), and hence, $q_m/w_m\to 1$, $m\to\infty$. Now, \eqref{e:ratio.Rqm.Rwm} implies $R_{q_m}/R_{w_m}\to1$ as $m\to\infty$. 
Recalling also $a\in \text{RV}_{1-\tau}$ and using the uniform convergence of regularly varying sequences, 
$$
\frac{a(R_{q_m})}{a(R_{w_m})} \sim \Big( \frac{R_{q_m}}{R_{w_m}} \Big)^{1-\tau} \to 1, \ \ \ m\to\infty; 
$$
hence, \eqref{e.ratio.aR.R} is obtained. 

Turning to \eqref{e:weak.core.cond.variant}, we note that by \eqref{e:weak.core.regime}, 
\begin{align*}
\binom{u_{m+1}}{2k+1} f(R_{q_m})^{2k+1} &\sim \frac{1}{(2k+1)!}\, \Big( \frac{u_{m+1}}{q_m} \Big)^{2k+1} \big( q_mf(R_{q_m}) \big)^{2k+1} \\
&\sim \frac{\xi^{2k+1}}{(2k+1)!}\, \Big( \frac{u_{m+1}}{q_m} \Big)^{2k+1} \to \frac{\xi^{2k+1}}{(2k+1)!}, \ \ \ m\to\infty, 
\end{align*}
where the last convergence is obtained as a result of \eqref{e:ratio.um.qm}. 

Returning to \eqref{e:ETm.normalized.light}, let us next calculate the limits for each term in the integral, while finding their appropriate upper bounds. Under our setup, it is straightforward to see that $a'\in \text{RV}_{-\tau}$ \citep[see, e.g., Proposition 2.5 in ][]{htresnick}. Therefore, $a(z)/z\to 0$ as $z\to\infty$, and for all $\rho>0$, 
\begin{equation}  \label{e:aRqm.Rqm.rho}
\Big( 1+\frac{a(R_{q_m})}{R_{q_m}}\rho \Big)^{d-1} \to 1, \ \  \ m\to\infty. 
\end{equation}
Note also that \eqref{e:aRqm.Rqm.rho} is bounded by $2(1\vee \rho)^{d-1}$ for sufficiently large $m$. 

Next we deal with $f\big( R_{q_m} + a(R_{q_m})\rho \big)/f(R_{q_m})$. Write 
\begin{align}
\frac{f\big( R_{q_m} + a(R_{q_m})\rho \big)}{f(R_{q_m})} &= \exp \Big\{  -\psi \big(  R_{q_m} + a(R_{q_m})\rho\big) + \psi(R_{q_m}) \Big\} \label{e:ratio.first.density} \\
&=\exp \bigg\{ -\int_0^\rho \frac{a(R_{q_m})}{a\big( R_{q_m} + a(R_{q_m})r \big)} \dif r  \bigg\}. \notag
\end{align}
By the uniform convergence of regularly varying functions and $a(z)/z\to0$ as $z\to\infty$, we have for every $r \geq 0$ that
$$
\frac{a(R_{q_m})}{a\big( R_{q_m} + a(R_{q_m})r \big)} \to 1, \ \ \ m\to\infty. 
$$
Therefore, for every $\rho> 0$, 
$$
\frac{f\big( R_{q_m} + a(R_{q_m})\rho \big)}{f(R_{q_m})} \to e^{-\rho},  \ \ \ m\to\infty. 
$$
Additionally, we define a sequence $\big( s_\ell(m), \, \ell \ge 0, \, m \ge 0 \big)$ by
$$
s_\ell(m) = \frac{\psi^\leftarrow \big( \psi(R_{q_m}) + \ell \big)-R_{q_m}}{a(R_{q_m})}, 
$$
equivalently, $\psi \big( R_{q_m} + a(R_{q_m})s_\ell(m) \big)=\psi(R_{q_m})+\ell$. Then, Lemma 5.2 in \cite{balkema2004} implies that for any $\epsilon\in (0,d^{-1})$, there exists a positive integer $N=N(\epsilon)$ such that $s_\ell(m) \le \epsilon^{-1}e^{\ell \epsilon}$ 
for all $m\ge N$ and $\ell \ge 0$. 
Since $\psi$ is increasing, we can establish the bound of \eqref{e:ratio.first.density} as follows: for $m\ge N$, 
\begin{align*}
&\exp \Big\{  -\psi\big( R_{q_m} + a(R_{q_m})\rho \big) + \psi(R_{q_m}) \Big\}\, \one \big\{ \rho>0 \big\}  \\
&=\sum_{\ell=0}^\infty \one \big\{ s_\ell(m) < \rho \le s_{\ell+1}(m) \big\} \exp \Big\{  -\psi\big( R_{q_m} + a(R_{q_m})\rho \big) + \psi(R_{q_m}) \Big\} \\
&\le \sum_{\ell=0}^\infty \one \big\{ 0 < \rho \le \epsilon^{-1} e^{(\ell+1)\epsilon} \big\} e^{-\ell}. 
\end{align*}

We now discuss the final untreated term from the integral in \eqref{e:ETm.normalized.light}. Let us give a helpful fact about $\big\| \big(R_{q_m} +a(R_{q_m})\rho  \big)\theta +y_i \big\|$ for $i\in \{ 1,\dots,2k \}$. We have that 
\begin{align*}
&\Big\| \big(R_{q_m} +a(R_{q_m})\rho  \big)\theta +y_i \Big\| - \Big(  R_{q_m} +a(R_{q_m})\rho +\langle \theta, y_i \rangle \Big) \\
&= \frac{\|y_i\|^2 - \langle \theta, y_i \rangle^2}{\Big\| \big(R_{q_m} +a(R_{q_m})\rho  \big)\theta +y_i \Big\|  + R_{q_m} +a(R_{q_m})\rho  + \langle \theta, y_i \rangle} =: \gamma_m(\rho, \theta, y_i). 
\end{align*}
In particular, if $\big\| \big(R_{q_m} +a(R_{q_m})\rho  \big)\theta +y_i \big\| \ge R_{q_m}$, then 
\begin{equation}  \label{e:unif.conv.gamma.m}
\big|  \gamma_m(\rho, \theta, y_i)\big| \le \frac{\big|\, \|y_i\|^2-\langle \theta, y_i \rangle^2  \, \big|}{2R_{q_m} + \langle \theta, y_i \rangle} \to 0, \ \ \ m\to\infty. 
\end{equation}
This convergence takes place uniformly for $\rho>0$, $\theta\in S^{d-1}$, and $y_i\in \R^d$ with $\|y_i\|\le ct$, where $c$ is determined by property \textbf{(H3)} of $h$---see Section \ref{s:prelim}. Continuing onward, let 
$$
A_m=\Big\{ y\in \R^d: \big\| \big( R_{q_m}  + a(R_{q_m})\rho \big)\theta +y \big\| \ge R_{q_m} \Big\};  
$$
then, for each $i\in \{1,\dots,2k\}$, 
\begin{align*}
&\frac{f\Big( \big\| \big( R_{q_m} +a(R_{q_m})\rho \big)\theta +y_i \big\| \Big)}{f(R_{q_m})}\, \one_{A_m}(y_i)  \\
&= \exp\Big\{ - \psi \big( R_{q_m} + a(R_{q_m})\rho + \langle \theta, y_i \rangle +\gamma_m(\rho, \theta, y_i) \big) + \psi(R_{q_m})\Big\}\, \one_{A_m}(y_i) \\
&= \exp\bigg\{  -\int_0^{\rho+\xi_m(\rho, \theta, y_i)} \frac{a(R_{q_m})}{a\big( R_{q_m} + a(R_{q_m})r \big)} \dif r   \bigg\}\, \one_{A_m}(y_i), 
\end{align*}
where $\xi_m(\rho,\theta, y_i) = a(R_{q_m})^{-1} \big( \langle \theta, y_i \rangle +\gamma_m(\rho,\theta,y_i) \big)$. 
Note that the last term is bounded by $1$, due to the fact that 
$$
\big\| \big( R_{q_m} +a(R_{q_m})\rho \big)\theta +y_i \big\| \ge R_{q_m} \ \Leftrightarrow \ \rho+\xi_m(\rho,\theta, y_i)\ge 0. 
$$
Additionally, \eqref{e:crackle.cond} and \eqref{e:unif.conv.gamma.m} yield that 
$$
\xi_m(\rho,\theta, y_i) \to \zeta^{-1} \langle \theta, y_i \rangle, \ \ \ m\to\infty, 
$$
for all $\rho>0$, $\theta\in S^{d-1}$, and $y_i\in \R^d$. Thus, as $m\to\infty$,  
$$
\exp\bigg\{  -\int_0^{\rho+\xi_m(\rho, \theta, y_i)} \frac{a(R_{q_m})}{a\big( R_{q_m} + a(R_{q_m})r \big)} \dif r   \bigg\} \to \exp \big\{ -\rho - \zeta^{-1} \langle \theta, y_i \rangle \big\}, 
$$
and 
$$
\one_{A_m}(y_i) \to \one \big\{ \rho + \zeta^{-1}\langle \theta, y_i \rangle \ge 0\big\}. 
$$

Combining all the bounds derived thus far, the integral in \eqref{e:ETm.normalized.light} is bounded above by 
\begin{align*}
&2\int_0^\infty \int_{S^{d-1}} \int_{(\R^d)^{2k}} h_t(0,\y)\, (1\vee \rho)^{d-1} \sum_{\ell=0}^\infty \one \big\{ 0<\rho \le \epsilon^{-1} e^{(\ell+1)\epsilon}\big\} e^{-\ell} \dif \y \dif \nu_{d-1}(\theta) \dif \rho \\
&=C^* \int_0^\infty \sum_{\ell=0}^\infty \one \big\{ 0<\rho \le \epsilon^{-1} e^{(\ell+1)\epsilon}\big\} e^{-\ell}  (1\vee \rho)^{d-1}\dif \rho \\
&\le C^* \Big( \frac{e^\epsilon}{\epsilon} \Big)^d \sum_{\ell=0}^\infty e^{-(1-\epsilon d)\ell} <\infty,
\end{align*}
as $\epsilon^{-1}e^{(\ell+1)\epsilon}\ge 1$ and $\epsilon d<1$. 
Now, by the dominated convergence theorem, we can see that the integral in \eqref{e:ETm.normalized.light} converges to 
\begin{align*}
\int_0^\infty \int_{S^{d-1}} \int_{(\R^d)^{2k}} h_t(0,\y)\, e^{-(2k+1)\rho - \zeta^{-1}\sum_{i=1}^{2k}\langle \theta, y_i \rangle } \prod_{i=1}^{2k}\one \big\{ \rho+\zeta^{-1} \langle \theta, y_i \rangle \ge 0 \big\} \dif \y \dif \nu_{d-1}(\theta) \dif \rho. 
\end{align*}
Because of this convergence, as well as \eqref{e.ratio.aR.R} and \eqref{e:weak.core.cond.variant}, we can get \eqref{e:first.Tm.light} as required. 

The proof of \eqref{e:first.Um.light} is almost the same as above, so we skip it. We can now conclude this lemma by showing \eqref{e:second.Tm.light} and \eqref{e:second.Um.light}. However, the proof has been omitted as it has essentially the same character as the proofs of \eqref{e:second.Tm} and \eqref{e:second.Um} in Lemma \ref{l:var_res}---albeit with different upper bounds as discussed above, due to the differing nature of the tail of probability densities. 
\end{proof}

\begin{proof}[Theorem \ref{t:lite}]
First, for every $n\ge1$, 
$$
\frac{U_m(t)}{a(R_{v_m})R_{v_m}^{d-1}} \le \frac{\chi_n^{(1)}(t)}{a(R_n)R_n^{d-1}} \le \frac{T_m(t)}{a(R_{w_m})R_{w_m}^{d-1}}. 
$$
Lemma \ref{l:var_res_lite} yields that 
$$
\limsup_{n\to\infty} \frac{\chi_n^{(1)}(t)}{a(R_n)R_n^{d-1}}  \le K_1(t) + \limsup_{m\to\infty} \frac{T_m(t)-\E[T_m(t)]}{a(R_{w_m})R_{w_m}^{d-1}}, 
$$
and 
$$
\liminf_{n\to\infty} \frac{\chi_n^{(1)}(t)}{a(R_n)R_n^{d-1}} \ge K_1(t) + \liminf_{m\to\infty}\frac{U_m(t)-\E[U_m(t)]}{a(R_{v_m})R_{v_m}^{d-1}}. 
$$
Now, the proof can be finished, if one can show that 
\begin{align}
&\big[ a(R_{w_m})R_{w_m}^{d-1} \big]^{-1} \big( T_m(t)-\E[T_m(t)] \big) \to 0, \ \ \ m\to\infty, \ \ \text{a.s.}, \label{e:conv.for.T_m.light}\\
&\big[ a(R_{v_m})R_{v_m}^{d-1} \big]^{-1} \big( U_m(t)-\E[U_m(t)] \big) \to 0, \ \ \ m\to\infty, \ \ \text{a.s.}  \label{e:conv.for.U_m.light}
\end{align}
By \eqref{e:second.Tm.light} in Lemma  \ref{l:var_res_lite}, for every $\epsilon>0$, 
\begin{align*}
\sum_{m=1}^\infty \P \Big( \, \big|\, T_m(t)-\E[T_m(t)] \, \big| > \epsilon a(R_{w_m})R_{w_m}^{d-1} \Big) \le \frac{1}{\epsilon^2}\sum_{m=1}^\infty \frac{\text{Var}\big( T_m(t) \big)}{\big(  a(R_{w_m})R_{w_m}^{d-1} \big)^2} \le C^* \sum_{m=1}^\infty \frac{1}{ a(R_{w_m})R_{w_m}^{d-1}}. 
\end{align*}
Because of the constraint in \eqref{e:constraint.gamma}, there exist $\delta_i>0$, $i=1,2$, so that 
$$
\gamma(d-\tau -\delta_1) \Big(\frac{1}{\tau} -\delta_2 \Big) > 1. 
$$
Then, $a\in \text{RV}_{1-\tau}$ implies that 
$$
a(R_{w_m}) R_{w_m}^{d-1} \ge C^* R_{w_m}^{d-\tau-\delta_1}
$$
for all $m\ge 1$. Note that by \eqref{e:R_n.asym}, 
\begin{align*}
R_{w_m} &\ge C^* \psi^\leftarrow \big( \log w_m \big) \ge C^* \psi^\leftarrow (\log u_m) \ge C^* m^{\gamma (1/\tau-\delta_2)}
\end{align*}
again for all $m\ge 1$. Therefore, 
$$
a(R_{w_m})R_{w_m}^{d-1} \ge C^* m^{\gamma (d-\tau -\delta_1)(1/\tau-\delta_2)}, 
$$
and 
$$
\sum_{m=1}^\infty \frac{1}{a(R_{w_m})R_{w_m}^{d-1}} \le C^*\sum_{m=1}^\infty \frac{1}{m^{\gamma (d-\tau -\delta_1)(1/\tau-\delta_2)}} < \infty. 
$$
Now, the Borel-Cantelli lemma completes the proof of \eqref{e:conv.for.T_m.light}. The proof of \eqref{e:conv.for.U_m.light} is the same by virtue of \eqref{e:second.Um.light} in Lemma  \ref{l:var_res_lite}. 

Before concluding the proof, we show finiteness of the limit in \eqref{e:fslln.exp}. Using property \textbf{(H3)} of $h$, 
\begin{align*}
\Big| \,  \sum_{k=0}^\infty (-1)^k s_k(t) \, \Big| &\le \sum_{k=0}^\infty \frac{s_{d-1} \xi^{k+1}}{(k+1)!} \int_0^\infty \int_{(\R^d)^k} \prod_{i=1}^k \one \big\{ \|y_i\| \le ct \big\}\, e^{-\rho} \dif \y \dif \rho \\
&\le C^* \sum_{k=0}^\infty \frac{\big( (ct)^d \xi \omega_d \big)^k}{k!} = C^* e^{(ct)^d \xi \omega_d } <\infty. 
\end{align*}
\end{proof}

\vspace{12pt}
\noindent \textbf{Acknowledgements}: The authors would like to thank the anonymous referee and the Associate Editor for their helpful insights, which have made the paper much more accessible. This research is partially supported by the National Science Foundation (NSF) grant, Division of Mathematical Sciences (DMS), \#1811428.

\bibliographystyle{abbrvnat}
\bibliography{EVT_EC_ProcessArXiv}


\end{document}